\documentclass{l4dc2027}
\title[SMC for Inequality-Constrained Optimization]{
Inequality-Constrained Optimization via Sliding Mode Control
}
\usepackage{times}
\usepackage{algorithm}
\usepackage{algpseudocode}
\usepackage{float}
\usepackage{enumitem}
\usepackage{comment}
\usepackage{listings}
\usepackage{verbatim}
\usepackage{bm}
\usepackage{graphicx}
\usepackage{subcaption}
\usepackage{wrapfig}

\usepackage[prependcaption,colorinlistoftodos]{todonotes}
\newtheorem{thm}{Theorem}[section]
\newtheorem{lem}{Lemma}[section]

\newtheorem{assumption}{Assumption}
\newtheorem*{assumption*}{Assumption}

\newcommand{\bd}{\begin{definition}} 
\newcommand{\ed}{\end{definition}}
\newcommand{\bt}{\begin{thm}} 
\newcommand{\et}{\end{thm}}
\let\eps\varepsilon

\newcommand{\R}{\mathbb{R}}

\newcommand{\bi}{\begin{itemize}} 
\newcommand{\ei}{\end{itemize}} 
\newcommand{\bds}{\begin{description}} 
\newcommand{\eds}{\end{description}} 
\newcommand{\beq}{\begin{equation}} 
\newcommand{\eeq}{\end{equation}} 

\newcommand{\norm}[1]{\|#1\|}

\newcommand{\setdef}[2]{\{#1 \; | \; #2\}}

\newcommand{\map}[3]{#1\colon #2 \rightarrow #3}

\newcommand{\e}{\mathrm{e}}

\newcommand{\sign}[1]{\operatorname{sign}\bigl({#1}\bigr)}

\newcommand{\1}{\mbox{\fontencoding{U}\fontfamily{bbold}\selectfont1}}
\newcommand{\0}{\mbox{\fontencoding{U}\fontfamily{bbold}\selectfont0}}

\newcommand{\cosmic}{COSMIC}

\definecolor{gnred}{RGB}{255,91,89}

\definecolor{gnred1}{RGB}{71,0,0} % 470000
\definecolor{gnred2}{RGB}{117,0,0} % 750000
\definecolor{gnred3}{RGB}{164,0,0} % a40000
\definecolor{gnred4}{RGB}{211,0,0} % d30000
\definecolor{gnred5}{RGB}{255,0,0} % FF0000
\definecolor{gnred6}{RGB}{255,42,34} % FF2a22
\definecolor{gnred7}{RGB}{255,91,89} % ff5b59 --- favorite

\definecolor{gnblue1}{RGB}{0,36,71}   % 002447
\definecolor{gnblue2}{RGB}{0,60,118}  % 003c76
\definecolor{gnblue3}{RGB}{0,85,164}
\definecolor{gnblue4}{RGB}{0,108,212}
\definecolor{gnblue4}{RGB}{0,108,212}
\definecolor{gnblue5}{RGB}{0,133,255}  % 0085ff
\definecolor{gnblue6}{RGB}{35,156,255} % 239cff
\definecolor{gnblue7}{RGB}{88,177,255} % 58b1ff

\definecolor{gnbrown1}{RGB}{71,27,0}  % 471b00
\definecolor{gnbrown2}{RGB}{117,45,0} % 752d00
\definecolor{gnbrown3}{RGB}{164,62,0} % a43e00
\definecolor{gnbrown4}{RGB}{211,80,0} % d35000
\definecolor{gnbrown5}{RGB}{255,97,0} % ff6100
\definecolor{gnbrown6}{RGB}{255,127,26} % ff7f1a
\definecolor{gnbrown7}{RGB}{255,155,86} % ff9b56

\RestyleAlgo{ruled}                    % rules above/below, caption on top
\SetAlgoCaptionSeparator{}             % "Algorithm 1 caption" instead of "Algorithm 1: caption"
\LinesNumbered
\SetNlSty{textrm}{}{:}                 % line numbers as "1:", "2:", ...
\SetAlgoNlRelativeSize{-1}             % small line numbers
\DontPrintSemicolon                    % don't print ";" at end of lines
\SetAlgoNoLine                         % no vertical bars
\SetFuncSty{textsc}                    % function names in small caps
\SetKwFor{For}{for}{do}{end for}
\SetKwIF{If}{ElseIf}{Else}{if}{then}{else if}{else}{end if}
\SetKwProg{Fn}{function}{}{end function}
\SetKwFunction{Step}{Step}
\SetKw{Ret}{return}
\SetKw{Brk}{break}

\author{%
 \Name{Rawan Hoteit} \Email{rhoteit@ethz.ch}
\\
\addr Automatic Control Laboratory (IfA), ETH Zürich, 8092 Zürich, Switzerland
 \AND
 \Name{Veronica Centorrino} \Email{vcentorrino@ethz.ch}\\
 \addr Automatic Control Laboratory (IfA), ETH Zürich, 8092 Zürich, Switzerland
 \AND
 \Name{Efe C. Balta} \Email{efe.balta@inspire.ch}\\
 \addr Inspire AG, Zürich and Automatic Control Laboratory (IfA), ETH Zürich, 8092 Zürich, Switzerland
 \AND
 \Name{John Lygeros} \Email{lygeros@ethz.ch}\\
 \addr Automatic Control Laboratory (IfA), ETH Zürich, 8092 Zürich, Switzerland
}
\begin{document}

\maketitle
\begin{abstract}%
We present \cosmic~(Constrained Optimization via Sliding Mode for Inequality Constraints), a continuous-time solver for inequality-constrained optimization. \cosmic~couples primal gradient-flow dynamics with multiplier dynamics derived from sliding mode control laws on a manifold defined by the smoothed Fischer--Burmeister function. 
We prove that \cosmic's trajectories reach this manifold in finite time from any primal initial point, after which the dynamics reduces to the gradient flow of a logarithmic barrier function. 
This yields a control-theoretic interpretation of barrier methods while allowing infeasible initialization.
For strongly convex problems, we further establish exponential convergence to the minimizer. To address the numerical stiffness that arises in discretization, we develop a modified discrete-time algorithm. We validate \cosmic~across various machine learning applications, demonstrating its ability to solve optimization problems accurately while scaling effectively to high-dimensional settings.
\end{abstract}


%\item Limit the main text (not counting acknowledgments or references) to 10 PMLR-formatted pages, using this template.
%\item Include {\em in the main text} enough details, including proof details, to convince the reviewers of the contribution, novelty and significance of the submissions.
%\item Make sure to cite references for related literature and for any claims that are not established in the paper using the \texttt{\textbackslash citep} and \texttt{\textbackslash citet} commands, e.g. \citep{levine2016end}, \citet{doyle2003guaranteed}, \citet[Theorem 2.2.1]{vershynin2018high}. We have provided a sample \texttt{.bib} file.

\section{Introduction}
We consider optimization problems of the form 
\begin{equation}
\begin{aligned}
\min_{x \in \R^n} ~& f(x)\\
\text{s.t.} ~& g(x) \le \0_m,
\end{aligned}
\label{eq:OP_ineq}
\end{equation}
where $\map{f}{\R^n}{\R}$ and $\map{g}{\R^n}{\R^m}$ are assumed to be continuously
differentiable throughout the paper. Inequality-constrained optimization problems, which subsume equality constraints, are prominent in various control, robotics, and machine learning applications, where the inequality constraints typically encode safety requirements or system and resource limitations.
% Literature Review
Classical approaches to solving~\eqref{eq:OP_ineq} include projection-based methods, primal-dual algorithms, and penalty methods.
The system-theoretic view of algorithms offers an alternative framework for analyzing and designing optimization algorithms through the lens of control theory~\citep{FD-ZH-GB-SB-JL-MM:24}.
This perspective has recently provided a control-theoretic interpretation of saddle-point flows, where the primal variables act as the plant of a dynamical system, and the Lagrange multipliers act as PID controllers that enforce constraint satisfaction~\citep{VC-SMF-SP-DR:25,VC-RH-EB-JL:2026}.
While PID controllers offer asymptotic convergence in these methods, sliding mode control (SMC) offers finite-time convergence and robustness to bounded uncertainty~\citep{HK:02}. Recent work applies SMC for equality-constrained problems~\citep{BD-SP-SK-TD:2026} with finite-time convergence. However, to date there is no first-order SMC-based solver
for inequality-constrained problems, with convergence guarantees.

To fill this gap, we propose an SMC-based solver for~\eqref{eq:OP_ineq}. The resulting continuous-time optimization algorithm is computationally efficient and enjoys finite-time convergence guarantees. Moreover, it provides a control-theoretic interpretation of penalty methods, which, to the best of our knowledge, has not been established in the literature.

\paragraph{Related Work.}
There has been extensive work in viewing optimization algorithms as continuous time closed-loop systems with feedback controllers~\citep{AH-ZH-SB-GH-FD:24, AA-JC:2023, GB-JC-JP-ED:2021}.
This perspective builds on the classical primal--dual (Arrow--Hurwicz--Uzawa) gradient flows~\citep{KJA-LH-HU:58}, in which the multipliers act as an integral controller on the constraint violation.~\cite{FC-WR:2020} utilize nonlinear sign-based control to design projected gradient flows for linear-equality constrained convex problems.
Within this line of work, \citet{VC-SMF-SP-DR:25, VC-RH-EB-JL:2026, CDC:VC-RH-EB-JL:2026} set the multipliers to be PID controllers on the primal flow.
Moreover, \citet{RZ-AR-JS-NL:25} apply feedback linearization to inequality-constrained problems and relate the resulting algorithm to sequential quadratic programming.  Closest to our work,~\cite{BD-SP-SK-TD:2026} designs the dual dynamics of equality constrained problems as SMC controllers and \citet{AZ-KY:2026} use discrete-time SMC to obtain a second-order solver for equality and inequality constrained problems, choosing a sliding manifold that encodes the KKT conditions and using the Fischer--Burmeister (FB) function.

\paragraph{Contributions.} 
We introduce the Constrained Optimization via Sliding Mode for Inequality Constraints (\cosmic), a continuous-time flow for solving~\eqref{eq:OP_ineq}. Under \cosmic, the primal variables follow a gradient descent on the Lagrangian, while the multipliers are driven by SMC laws whose sliding manifold is the zero set of the smoothed FB function of the constraints and multipliers. On this manifold, primal and dual feasibility hold strictly and complementary slackness holds up to a perturbation $\epsilon>0$. At equilibrium, the primal dynamics enforce stationarity, and the equilibria satisfy a perturbed KKT system, which reduces to the exact one as $\epsilon\to0$.
%converge guarantees

Our main results establish finite-time convergence to the sliding manifold and, for strongly convex problems, exponential convergence to the minimizer of~\eqref{eq:OP_ineq}.
We also show that, once on the manifold, \cosmic~reduces to the gradient flow of a logarithmic barrier method.
This yields a control-theoretic interpretation of barrier methods and generalizes them, in the sense that our solver does not require a feasible
initialization.
To address the numerical difficulties that arise when discretizing \cosmic, we propose a modified discrete-time implementation.
%numerics
Finally, through training a support vector machine (SVM) on the Iris Setosa-vs-Versicolor dataset, and image completion with an input convex neural network (ICNN), we  highlight the computational efficiency and accuracy of \cosmic~relative to other nonlinear control-based optimization solvers. The code to reproduce the experiments is available at~\url{https://github.com/sudo-Ra-ONE/COSMIC}.

\begin{figure}
    \centering
\includegraphics[width=0.9\linewidth]{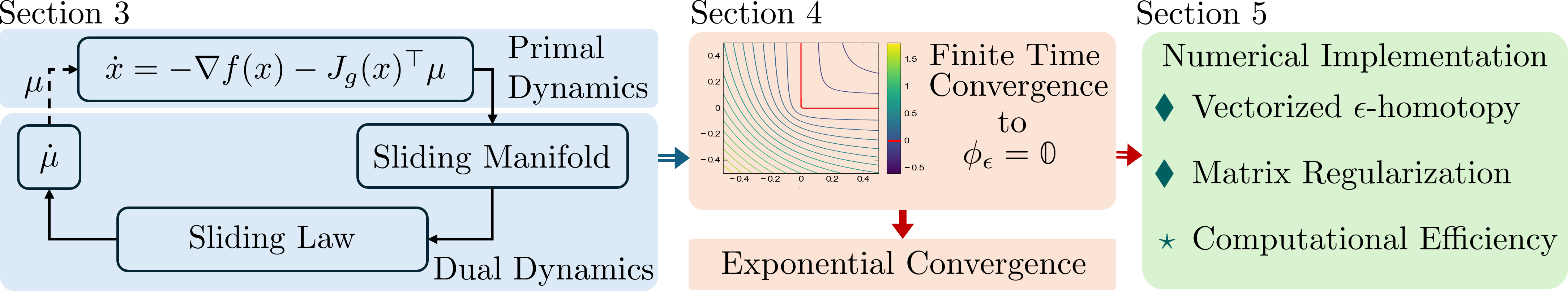}
    \caption{Visual Outline of the Paper}
    \label{fig: Visual Outline}
\end{figure}

\paragraph{Outline.} Section~\ref{sec: preliminaries} presents the mathematical preliminaries and is followed by the problem formulation and design of
%the constrained optimization solver, 
COSMIC, in Section~\ref{sec: problem_setup}. Section~\ref{sec: results} formalizes the analysis and the numerical implementation of our solver, and Section~\ref{sec: numerical_results} benchmarks it on various problems. Section~\ref{sec: conclusion} provides a summary and directions for future work. We provide a visual outline in Figure~\ref{fig: Visual Outline}.

 \label{sec: introduction}
\section{Mathematical Preliminaries}
\paragraph{Notation.} We denote by $\R$ and $\R_{>0}$ the real and positive real numbers, respectively, by $\norm{\cdot}$ the Euclidean norm on $\R^n$, and by $\0_n$ and $\1_n$ the all-zeros and all-ones vectors in $\R^n$.
%The sets of real and positive real numbers are denoted by $\R$ and $\R_{>0}$, respectively. For $n \in \mathbb{N}$, $\R^n$ denotes the $n$-dimensional Euclidean space, equipped with the Euclidean norm $\norm{\cdot}$. We let $\0_n$ and $\1_n$ be the all-zeros and all-ones vectors in $\R^n$, respectively. 
We let $\oslash$ be the per-component division between two vectors.
%Given $A, B \in \R^{n\times n}$ symmetric, we write $A \preceq B$ (resp. $A \prec B$) if $B-A$ is positive semidefinite (resp. definite). For a symmetric matrix $A$, we let $\subscr{\lambda}{min}(A)$ and $\subscr{\lambda}{max}(A)$ be its minimum and maximum eigenvalue, respectively.
The class of $k$-times continuously differentiable functions is $\mathcal{C}^k$. For a differentiable $\map{f}{\R^n}{\R}$, $\nabla f(x)$ is its gradient.
%Let $\map{f}{\R^n}{\R}$ be twice differentiable. Then $\nabla f(x)$ and $\nabla^2 f(x)$ denote its gradient and Hessian matrices, respectively.
%For $\map{f}{\R^{n+m}}{\R}$, we let $\nabla_{x} f(x,y) := \frac{\partial f}{\partial x}(x,y)$ be its partial gradients with respect to $x$.
Given $\map{h}{\R^n}{\R^m}$ continuously differentiable, $J_h(x) \in \R^{m \times n}$ denotes its Jacobian matrix.
Given $\rho, L > 0$, a map $\map{f}{\R^n}{\R}$ is (i) \emph{$\rho$-strongly convex} if $x \mapsto f(x) - \frac{\rho}{2}\norm{x}^2$ is convex; (ii) \emph{$L$-smooth} if it is differentiable and $\nabla f$ is $L$-Lipschitz. The sign function is defined by $\sign{x} := 1$ if $x>0$, $\sign{x}:=-1$ if $x<0$, and $\sign{0}:=0$, and it is applied component-wise to vectors.
\paragraph{Optimality Conditions.}
 Let $\mathcal{F} := \setdef{x \in \R^n}{g(x) \leq \0_m}$ denote the feasible set of~\eqref{eq:OP_ineq}, $f^\star := \inf_{x\in\mathcal F} f(x)$ its optimal value, with the convention
$\inf\emptyset=+\infty$, and $\mathcal{A}(x) := \setdef{i \in \{1,\dots,m\}}{g_i(x) = 0}$ denote the active set at $x \in \mathcal{F}$. 
We say that linear independence constraint qualification (LICQ) holds at $x \in \mathcal{F}$ if the gradients $\{\nabla g_i(x)\}_{i \in \mathcal{A}(x)}$ are linearly independent. 
A pair $(x^\star, \mu^\star) \in \R^{n+m}$ satisfies the Karush--Kuhn--Tucker (KKT) conditions if
\begin{subequations}
\label{eq:kkt_conditions}
\begin{align}
&\nabla f(x^\star) + J_g(x^\star)^\top \mu^\star = \0_n,
\label{eq:kkt_stationarity}\\
&g(x^\star) \le \0_m, \quad \mu^\star \ge \0_m, \quad \mu^\star_i g_i(x^\star) = 0,\quad \forall i=1,\dots,m.
\label{eq:kkt_complementarity}
\end{align}
\end{subequations}

For $\epsilon \geq 0$, let $\map{\phi_\epsilon}{\R \times \R}{\R}$, $\phi_\epsilon (a,b) := a+b-\sqrt{a^2+b^2+\epsilon^2}$. The Fischer--Burmeister (FB) function~\citep{AF:92} is $\phi_0$, and $\phi_\epsilon$, $\epsilon>0$, is its smoothed version~\citep{CK:96}. For $\epsilon>0$, $\phi_\epsilon\in\mathcal C^1$ and $\phi_\epsilon(a,b) = 0$ if and only if $a, b>0$ and $ab = \epsilon^2/2$. Hence, with $a := -g_i(x)$ and $b := \mu_i$, $\phi_\epsilon$ encodes a perturbed version of the conditions in \eqref{eq:kkt_complementarity}, which are exact for $\epsilon=0$.
%For $\epsilon \geq 0$, the smoothed Fischer--Burmeister (FB) function~\citep{to_add} is $\map{\phi_\epsilon}{\R \times \R}{\R}$, $\phi_\epsilon (a,b) := a+b-\sqrt{a^2+b^2+\epsilon^2}$, and the standard FB function is $\phi_0$. For $\epsilon>0$, $\phi_\epsilon\in\mathcal C^1$ and $\phi_\epsilon(a,b) = 0$ if and only if $a>0$, $b>0$, and $ab = \epsilon^2/2$. Hence, applied to $a := - g_i(x)$ and $b = \mu_i$, $\phi_\epsilon$ encodes a perturbed version of~\eqref{eq:kkt_primal_feasibility}--\eqref{eq:kkt_complementarity}, which reduces to the exact conditions when $\epsilon=0$.
%Throughout this work, the sign function is defined by $\sign{0}=0$.
% \paragraph{Sliding Mode Control.}
% Sliding mode control (SMC) is a nonlinear control technique that drives a system's state trajectories into a predefined invariant surface called the sliding manifold~\citep{HK:02}. The design of a sliding mode controller is based on two phases: the reaching phase that drives the states into the sliding manifold as dictated by a sliding law, and the sliding phase in which the states remain on the manifold and slide towards the equilibrium point. 
 \label{sec: preliminaries}
\section{Problem Formulation}
We treat the primal dynamics of~\eqref{eq:OP_ineq} as the plant and design the Lagrange multipliers as SMC-based feedback controllers that enforce the KKT conditions.
SMC drives the system states to a predefined invariant, sliding manifold~\citep{HK:02}, in two phases: a reaching phase, and a sliding phase, during which they slide toward the equilibrium.
We next make this construction precise. 

Consider the Lagrangian associated with~\eqref{eq:OP_ineq}, that is, the map $\map{L}{\R^n \times \R^m}{\R}$ defined by
$L(x, \mu)=f(x) + \mu^\top g(x),$ where $\mu \in \R_{\geq 0}^m$ is the Lagrange multiplier of the constraints.
Treating $\mu$ as a control input to the gradient flow of $L$ with respect to $x$ yields the system 
\beq
\label{eq:open_system_ineq}
\dot{x}= - \nabla_x L(x, \mu) = - \nabla f(x) - J_g(x)^\top \mu,
\eeq
%We achieve the KKT conditions, through setting $\mu\in\R^m$ as the SMC input to the primal dynamics. 
with output $y = g(x)$. We design the input $\mu$ via SMC so that the trajectories of~\eqref{eq:open_system_ineq} reach, and remain on, a manifold on which an $\epsilon$-perturbed KKT conditions hold.
Fix $\epsilon>0$. With a slight abuse of notation, we consider $\map{\phi_\epsilon}{\R^n\times\R^m}{\R^m}$, defined by $\phi_\epsilon(x,\mu) := \big(\phi_\epsilon(-g_1(x),\mu_1),\dots,\phi_\epsilon(-g_m(x),\mu_m)\big)$.
Let the sliding manifold be given by $S(x,\mu) = \setdef{(x,\mu) \in \R^{n+m}}{\phi_\epsilon(x,\mu) = \0_m}$, where $S_i = \phi_{\epsilon, i}(x,\mu)$ for $i\in\{1,\dots,m\}$, and consider the following sliding law
\beq \label{eq:Sliding law}
\dot{S}_i = -KS_i-\beta\sign{S_i},
\eeq 
where $K >0$ and $\beta \geq 0$. Differentiating $S$ yields
\beq\label{eq:S_dot}
\dot S = \dot \phi_\epsilon = D_\epsilon(x,\mu) \dot \mu-G_\epsilon (x,\mu)J_g\dot x,
\eeq
where $D_\epsilon(x,\mu), G_\epsilon(x,\mu) \in \R^{m\times m}$ are diagonal with entries $D_i = 1-\frac{\mu_i}{\sqrt{g_i(x)^2+\mu_i^2+\epsilon^2}}$ and $G_i = 1+\frac{g_i(x)}{\sqrt{g_i(x)^2+\mu_i^2+\epsilon^2}}$.
%Since $\epsilon>0$, $|\mu_i|<\sqrt{g_i(x)^2+\mu_i^2+\epsilon^2}$ and $|g_i(x)| < \sqrt{g_i(x)^2+\mu_i^2+\epsilon^2}$, hence 
Since $D_i, G_i \in ]0,2[$ for all $i$ and all $(x,\mu)$, $D_\epsilon(x,\mu)$ is invertible.
Combining~\eqref{eq:open_system_ineq} and~\eqref{eq:S_dot} yields the closed loop dynamics  
\begin{equation}
\label{eq:SMC-Dyn}
\begin{cases}
 \dot x = -\nabla f(x)-J_g^\top(x) \mu,\\
 \dot{\mu} =D_{\epsilon}^{-1}(x, \mu) \left(-K \phi_\epsilon(x, \mu) - \beta \sign{\phi_\epsilon(x, \mu)}+ G_{\epsilon}(x, \mu) J_g(x) \dot{x}\right),
\end{cases}
\end{equation}
which we refer to as \emph{\cosmic}.
Because~\eqref{eq:SMC-Dyn} is discontinuous, we consider solutions of~\eqref{eq:SMC-Dyn} in the sense of Filippov~\citep{JC:2008}. 
Since $f,g\in\mathcal C^1$ and $D_\epsilon^{-1}$ is continuous, the right-hand side of~\eqref{eq:SMC-Dyn} is measurable and locally essentially bounded.
%that is, bounded on a bounded neighborhood of every point, excluding sets of measure zero.
Then, for every $(x_0,\mu_0)\in\R^{n+m}$, a Filippov solution of~\eqref{eq:SMC-Dyn} exists~\citep[Proposition~3]{JC:2008} on a maximal interval $[0,T_{\max})$, with $T_{\max}\in\,]0,\infty]$.

    % Proving that the xiltiplier remains non-negative
    % Considering the closed loop system in~\eqref{eq:OP_ineq} and let $x(t)$ be a trajectory of~\eqref{eq:open_system_ineq}. If $\mu(0) \in \R^m_{\geq 0}$, then $\mu(t) \in \R^m_{\geq 0}$ for all $t \geq 0$. This can be seen through considering $\dot\mu_i = {F_{\text{r},i}}$, where $F_{\text{r},i} = \relu{M_i^{-1} \left( -K \phi_{\epsilon,i} - \beta \sign{\phi_{\epsilon,i}}+  G J_g \dot{x} \right)}$. For all $\mu \in \R_{\geq 0}^m$ such that $\mu_i =0$ we have $F_{\text{r},i}(x,\mu) \geq 0.$ The statement then follows from Nagumo's Theorem~\cite{MN:1942}. For $\mu \geq 0$, and for all $x\in \mathbb{R}^m$, $\phi_\epsilon \leq 0$.

\bt[Equilibria of~\eqref{eq:SMC-Dyn} and KKT points of~\eqref{eq:OP_ineq}]
\label{thm:eq_KKT}
%Assume that $f$, $g$ are continuously differentiable, and c
Consider the \cosmic~dynamics~\eqref{eq:SMC-Dyn} with $K, \eps >0$, and $\beta \geq 0$. Let $\Omega_\epsilon$ denote the set of equilibria of~\eqref{eq:SMC-Dyn}. The following statements hold.
\begin{enumerate}[label=\textup{(\roman*)}]
\item \label{th_eqKKT:item1} %check also prop 4.1 in Bertsekas, Nonlinear programming for barrier methods.
\emph{(Perturbed KKT characterization)}
For every $\epsilon>0$, $(\hat{x},\hat{\mu}) \in \Omega_\epsilon$ if and only if $\nabla f(\hat{x})+ J_g(\hat{x})^\top \hat \mu=\0_n$ and $\phi_\epsilon(\hat{x},\hat{\mu})=\0_m$.
\item \label{th_eqKKT:item2}
\emph{(Limiting KKT points)}
Let $(x_k, \mu_k) \in \Omega_{\epsilon_k}$ and $\epsilon_k \to 0$. Then, every limit point $(\bar x,\bar\mu)$ of $\{(x_k,\mu_k)\}$ satisfies the KKT conditions of~\eqref{eq:OP_ineq}.
\item \label{th_eqKKT:item3}
\emph{(Multiplier convergence under LICQ)}
Let $(x_k,\mu_k)\in \Omega_{\epsilon_k}$, $\epsilon_k \to 0$, and let $\bar x$ be a limit point of $\{x_k\}$ at which LICQ holds. Then, along every subsequence such that $x_k \to \bar x$, $\{\mu_k\}$ converges to the unique $\bar\mu$ such that $(\bar x,\bar\mu)$ satisfies the KKT conditions of~\eqref{eq:OP_ineq}.
\item \label{th_eqKKT:item4}
\emph{(Boundedness and suboptimality in the convex case)}
If $f$ and $g_i$, $i\in\{1,\dots,m\}$, are convex and the set of minimizers of~\eqref{eq:OP_ineq} is nonempty and compact, then, for every $\bar\epsilon>0$, the set of all equilibria obtained for $\epsilon\in\,]0,\bar\epsilon]$, i.e., $\bigcup_{\epsilon\in]0,\bar\epsilon]} \Omega_\epsilon$, is bounded, and
\begin{equation}\label{eq:subopt}
0 \leq f(x)-f^\star \leq \dfrac{m\epsilon^2}{2} \qquad \text{for all }(x,\mu) \in \Omega_\epsilon .
\end{equation}
Consequently, any sequence $\{(x_k,\mu_k)\}$ as in~\ref{th_eqKKT:item2} is bounded and has at least one limit point satisfying the KKT conditions of~\eqref{eq:OP_ineq}.
\end{enumerate}
\et
\begin{proof}
%i
\emph{Item~\ref{th_eqKKT:item1}.}
Fix $\epsilon>0$. Since $D_\epsilon(x,\mu)$ is invertible, $(\hat x,\hat\mu)\in\Omega_\epsilon$ if and only if $\nabla f(\hat x)+J_g(\hat x)^\top\hat\mu=\0_n$ and, using $\dot{\hat x}=\0_n$, $K\phi_\epsilon(\hat x,\hat\mu)+\beta\sign{\phi_\epsilon(\hat x,\hat\mu)}=\0_m$. Since $K>0$ and $\beta\geq 0$, by definition of $\sign{\cdot}$, the last equality holds if and only if $\phi_\epsilon(\hat x,\hat\mu)=\0_m$, which in turn is equivalent to $g(\hat x)<\0_m$, $\hat\mu>\0_m$, and $\hat\mu_i g_i(\hat x)=-\epsilon^2/2$ for all $i$.

% ii
\emph{Item~\ref{th_eqKKT:item2}.} Let $(\bar x,\bar\mu)$ be a limit point of $\{(x_k,\mu_k)\}$. %In the remaining parts of the proof, 
With a slight abuse of notation, we denote by $\{(x_k,\mu_k)\}$ also a subsequence converging to $(\bar x,\bar\mu)$. Since $(x_k,\mu_k) \in \Omega_{\epsilon_k}$, item~\ref{th_eqKKT:item1} implies $\nabla f(x_k) + J_g(x_k)^\top\mu_k = \0_n$, $\mu_{k,i}g_i(x_k) = -\epsilon_k^2/2$, $g(x_k) < \0_m$, and $\mu_k > \0_m$, for all $i$ and $k$. By continuity of $\nabla f$, $J_g$, $g$, passing to the limit yields~\eqref{eq:kkt_conditions} at $(\bar x,\bar\mu)$.
%as $k \to \infty$ we have $\nabla f(\bar x)+J_g(\bar x)^\top\bar \mu = \0_n$, $g(\bar x)  \leq \0_m$, $ \bar \mu \geq \0_m$ and $\bar\mu_i g_i(\bar x)=0$ for all $i$, i.e., the KKT conditions at $(\bar x,\bar\mu)$. This shows statement~\ref{th_eqKKT:item2}.

% iii
\emph{Item~\ref{th_eqKKT:item3}.} Let $(x_k,\mu_k)\in \Omega_{\epsilon_k}$, $\epsilon_k \to 0$, and let $\bar x$ be a limit point of $\{x_k\}$ at which LICQ holds.
Item~\ref{th_eqKKT:item1} and continuity of $g$ implies %$g(x_k)<\0_m$, and hence 
$g(\bar x)\le\0_m$. Define the inactive constraints $\mathcal I :=\{1,\dots,m\}\setminus\mathcal{A}$, and let $J_{\mathcal A}(x)$ and $J_{\mathcal I}(x)$ be the rows of $J_g(x)$ indexed by $\mathcal{A}$ and $\mathcal{I}$, respectively, and use the same indexing also for $\mu$. By stationarity at $(x_k,\mu_k)$ we have
\beq
\label{eq:stat_split}
\nabla f(x_k)+J_{\mathcal A}(x_k)^\top\mu_{k,\mathcal A}+J_{\mathcal I}(x_k)^\top\mu_{k,\mathcal I} = \0_n.
\eeq
%By LICQ, $J_{\mathcal A}(\bar x)$ has full row rank, so $\tilde J(\bar x) := J_{\mathcal A}(\bar x)J_{\mathcal A}(\bar x)^\top$ is positive definite. By continuity of $x\mapsto \tilde J(x)$, $\tilde J(x_k)$ is positive definite, hence invertible, for $k$ large enough, and $\tilde J(x_k)^{-1} \to \tilde J(\bar x)^{-1}$. 
By LICQ, $J_{\mathcal A}(x)J_{\mathcal A}(x)^\top$ is positive definite at $\bar x$ and, by continuity, at some $x_k$ for $k$ large enough. Multiplying both sides of~\eqref{eq:stat_split} by $J_{\mathcal A}(x_k)$ yields
$
\mu_{k,\mathcal A} {=} - (J_{\mathcal A}(x_k)J_{\mathcal A}(x_k)^\top )^{-1} J_{\mathcal A}(x_k)\big(\nabla f(x_k)+J_{\mathcal I}(x_k)^\top\mu_{k,\mathcal I}\big).
$
Moreover, the inactive constraints' multipliers vanish in the limit as $k\to\infty$, i.e., as $x_k\to\bar x$ and $\epsilon_k\to 0$ jointly. To see this, let $i\in\mathcal I$, so that $g_i(\bar x)<0$. By continuity, $g_i(x_k)\leq g_i(\bar x)/2 < 0$, for $k$ large enough, and by item~\ref{th_eqKKT:item1}, $\mu_{k,i}g_i(x_k)=-\epsilon_k^2/2$ and $\mu_{k,i} > 0$. Then, we have
$ 0\!<\!\mu_{k,i}\!=\!\frac{\epsilon_k^2}{-2 g_i(x_k) }\!\leq\!\frac{\epsilon_k^2}{-g_i(\bar x)} \to 0 \text{ as } k \to \infty,$ resulting in $\bar \mu_I = 0$,
and, by continuity of $\nabla f$ and $J_g$,
$ \mu_{k,\mathcal A} \longrightarrow \bar \mu_{\mathcal A}:= - (J_{\mathcal A}(\bar x)J_{\mathcal A}(\bar x)^\top)^{-1}J_{\mathcal A}(\bar x)\nabla f(\bar x). $
Hence $\mu_k\to\bar\mu$, with $\bar \mu = (\bar \mu_A, \bar \mu)$. Finally, passing to the limit in~\eqref{eq:stat_split} yields $\nabla f(\bar x)+J_g(\bar x)^\top\bar\mu=\0_n$. Moreover, $\bar\mu\ge\0_m$ is a limit of positive vectors with $g(\bar x)\le\0_m$, and $\bar\mu_i g_i(\bar x)=0$ for all $i$, since $g_i(\bar x) = 0$ for $i \in \mathcal{A}$ and $\bar\mu_i =0$ for $i \in \mathcal{I}$. %This concludes the proof of item~\ref{th_eqKKT:item3}.

%iv
\emph{Item~\ref{th_eqKKT:item4}.} %Since the set of minimizers is nonempty, $f^\star$ is finite. 
Let $\epsilon>0$ and $(x,\mu)\in \Omega_\epsilon$. By item~\ref{th_eqKKT:item1}, $g(x)<\0_m$, $\mu > \0_m$, $\mu_i g_i(x) = -\epsilon^2/2$ for all $i$. Hence $x$ is feasible and $f(x) \geq f^\star$. Since $\mu>\0_m$ and $f$, $g_i$ are convex, $z\mapsto L(z,\mu)$ is convex and, by stationarity, $\nabla_x L(x,\mu) = \0_n$. Hence, $x$ minimizes $L(\cdot,\mu)$ and, for all $z\in\R^n$,
\begin{equation}\label{eq:4}
f(z) + \sum_{i=1}^m \mu_i g_i(z) = L(z,\mu) \geq L(x,\mu) = f(x) - \frac{m\epsilon^2}{2}.
\end{equation}
%For any $z\in\R^n$, we have
%\begin{align}
%f(z) &\geq f(x) + \nabla f(x)^\top(z - x) = f(x) - \sum_{i = 1}^m \mu_{i}\nabla g_i(x)^\top(z - x) \label{eq:2}\\ & \geq f(x) - \sum_{i = 1}^m \mu_{i} \big(g_i(z) - g_i(x)\big) = f(x) - \sum_{i = 1}^m \mu_{i} g_i(z) - \dfrac{m \epsilon^2}{2}, \label{eq:4}
%f(z) &\geq f(x) - \sum_{i = 1}^m \mu_{i}\nabla g_i(x)^\top(z - x) \geq f(x) - \sum_{i = 1}^m \mu_{i} g_i(z) - \dfrac{m \epsilon^2}{2}, \label{eq:4}\end{align}
%where we used convexity of $f$ and stationarity of $(x,\mu)$ first, convexity of $g_i$ and the fact that $\mu_{i} g_i(x)=-\epsilon^2/2$, for all $i$, then.
%Summing up, we have
%\beq
%\label{eq:gap}
%f(x) \leq f(z) + \sum_{i = 1}^m \mu_{k,i} g_i(z)+\dfrac{m\epsilon^2}{2}, \quad \text{ for all } z \in \R^n.
%\eeq
%Letting $z=x^\star$ be a minimizer, then $\mu_{i} g_i(x^\star) \leq 0$ for all $i$, so $f(x) - f^\star \leq \tfrac{m \epsilon^2}{2}$, which, together with $f(x)\geq f^\star$, proves~\eqref{eq:subopt}.
Choosing $z=x^\star$ as a minimizer, and using $\mu_i g_i(x^\star)\leq 0$ yields~\eqref{eq:subopt}.
Now let $\bar\epsilon>0$ and $(x,\mu)\in \Omega_\epsilon$ with $\epsilon\in]0,\bar\epsilon]$. Then $x \in \setdef{x \in \mathcal{F}}{f(x) \leq f^\star + m\bar\epsilon^2/2}$, which is bounded as a sublevel set of a convex problem whose set of minimizers is nonempty and compact~\citep[Corollary 8.7.1]{RTR:70}. Hence, $x$ is bounded uniformly in $\epsilon\in]0,\bar\epsilon]$. 
Moreover, fixing any $(x_\textup{s},\mu_\textup{s})\in\bigcup_{\epsilon\in]0,\bar\epsilon]}\Omega_\epsilon$, so that $g(x_\textup{s})<\0_m$ by item~\ref{th_eqKKT:item1}, and choosing $z = x_\textup{s}$, inequality~\eqref{eq:4} and $f(x) \geq f^\star$ yield
$
\sum_{i = 1}^m\mu_i \big(-g_i(x_\textup{s})\big) \leq f(x_\textup{s}) - f^\star + \tfrac{m\bar\epsilon^2}{2}.
$
Since $-g_i(x_\textup{s})>0$ and $\mu_i>0$ for all $i$, each term of the sum is nonnegative, so $0<\mu_i\le \frac{f(x_\textup{s}) - f^\star + m\bar\epsilon^2/2}{-g_i(x_\textup{s})}$, $i=1,\dots,m$, uniformly in $\epsilon \in ]0,\bar\epsilon]$. Therefore, $\bigcup_{\epsilon\in]0,\bar\epsilon]} \Omega_\epsilon$ is bounded.
Finally, let $(x_k,\mu_k)\in \Omega_{\epsilon_k}$ with $\epsilon_k\to 0$, and let $\bar\epsilon:=\sup_k\epsilon_k$, which is finite since $\{\epsilon_k\}$ converges. Then $\epsilon_k\in\,]0,\bar\epsilon]$ for all $k$, so $\{(x_k,\mu_k)\}\subset\bigcup_{\epsilon\in]0,\bar\epsilon]}\Omega_\epsilon$, which is bounded by the argument above. Then, by the Bolzano--Weierstrass theorem, $\{(x_k,\mu_k)\}$ admits a limit point, which satisfies the KKT conditions of~\eqref{eq:OP_ineq} by item~\ref{th_eqKKT:item2}. 
\end{proof}
Theorem~\ref{thm:eq_KKT}.\ref{th_eqKKT:item1} implies that $\Omega_\epsilon$ does not depend on the sliding mode parameters $K>0$ and $\beta \geq0$.
%Next, we note that statements~\ref{th_eqKKT:item2} and~\ref{th_eqKKT:item3} differ in what is assumed about the multipliers.
{Moreover, items~\ref{th_eqKKT:item1}--\ref{th_eqKKT:item3} assume neither feasibility of~\eqref{eq:OP_ineq} nor existence of a minimizer, they characterize the equilibria of~\eqref{eq:SMC-Dyn} whenever these exist. Item~\ref{th_eqKKT:item1} yields $g(x)<\0_m$ for every $(x,\mu)\in\Omega_\epsilon$, so $\Omega_\epsilon=\emptyset$ if~\eqref{eq:OP_ineq} has no strictly
feasible point, e.g., if it is infeasible. If $f$ and $g_i$ are convex, inequality~\eqref{eq:4}, which only uses convexity, implies that $\Omega_\epsilon=\emptyset$ if $f^\star=-\infty$.}
Statement~\ref{th_eqKKT:item2} requires no constraint qualification, since it assumes that $\{(x_k,\mu_k)\}$ has a limit point, which in particular requires $\{\mu_k\}$ to have a bounded subsequence.
On the other hand, item~\ref{th_eqKKT:item3} requires only that $\{x_k\}$ has a limit point and shows that LICQ at $\bar x$ guarantees convergence of $\{\mu_k\}$ to the (unique) KKT multiplier at $\bar x$. Without constraint qualification, the sequence $\{\mu_k\}$ may diverge even if $\{x_k\}$ converges to a minimizer. Consider, for example, $\min_{x\in\R}\{x : -x^3 \leq 0\}$, whose minimizer $x^\star=0$ violates LICQ, since $\nabla g(0)=0$, and admits no KKT multiplier. The equilibrium of~\eqref{eq:SMC-Dyn} is $(x_\epsilon,\mu_\epsilon) = \big(3\epsilon^2/2, 4/(27\epsilon^4)\big)$, so $x_\epsilon\to x^\star$ while $\mu_\epsilon\to\infty$ as $\epsilon \to 0$.
Finally, item~\ref{th_eqKKT:item4} specializes the results to convex problems and states that every equilibrium of~\eqref{eq:SMC-Dyn} is strictly feasible and $m\epsilon^2/2$-suboptimal for any fixed $\epsilon>0$, rather than only in the limit $\epsilon\to 0$. In particular, choosing $\epsilon \leq \sqrt{2\delta/m}$ guarantees $f(x)-f^\star \leq \delta$ for a prescribed accuracy $\delta>0$. Moreover, under the assumptions of item~\ref{th_eqKKT:item4} the boundedness of $\{\mu_k\}$ implicitly assumed in item~\ref{th_eqKKT:item2} holds automatically and the existence of a limit point satisfying the KKT conditions is guaranteed.% without requiring LICQ.
 \label{sec: problem_setup}
\section{Convergence Analysis}

%We study the convergence properties of \cosmic. We 
First, we show that trajectories of \cosmic~reach the sliding manifold $\mathcal{S}$ in finite time.
\begin{lem}[%Finite time convergence of~\eqref{eq:SMC-Dyn}
Finite-time reaching of the sliding manifold]
\label{lem:finite_time}
Let $\beta > 0$, and let $(x(t), \mu(t))$ be a Filippov solution of~\eqref{eq:SMC-Dyn} on $[0,T_{\max})$. Define $\phi_{i,0} := \phi_{\epsilon,i}(x(0),\mu(0))$, for $i\in\{1,\dots,m\}$, let
$t_\textup{r}  := \max_i\frac{1}{K}\ln\Bigl(\frac{K|\phi_{i,0}|+\beta}{\beta}\Bigr)$ and assume that $T_{\max}\!>\!t_\textup{r}$. Then $\phi_\epsilon(x(t), \mu(t))\!=\!\0_m$ for all $t\in[t_\textup{r},T_{\max})$.
%Under our choice of sliding law~\eqref{eq:Sliding law}, $\phi_\epsilon(x(t), \mu(t)) = \0_m$ for all $t\geq t_\textup{r}$ where $t_\textup{r}  = \max_i\frac{1}{K}\ln(\frac{K|\phi_{i,0}|+\beta}{\beta})$ and $\phi_\epsilon(x(t), \mu(t)) = \0_m$ is invariant thereafter.
\end{lem}
\begin{proof}
Since $\phi_{\epsilon,i}(x(t),\mu(t))\in\mathcal{C}^1$ and $(x(t),\mu(t))$ is absolutely continuous, $\phi_i(t) := \phi_{\epsilon,i}(x(t),\mu(t))$ is absolutely continuous and $\dot\phi_i = D_{\epsilon,i}\dot\mu_i - G_{\epsilon,i}\nabla g_i(x)^\top\dot x$ satisfies the sliding law~\eqref{eq:Sliding law} with $\dot\mu$ in~\eqref{eq:SMC-Dyn} that is, we have  $\dot\phi_i = -K\phi_i - \beta \sign{\phi_i}$, for almost every $t\in[0,T_{\max})$. In particular, this implies that each $\phi_i$ evolves independently of $x$ and of $\phi_j$, for every $j\neq i$.
If $\phi_{i,0}>0$, then, as long as $\phi_i(t)>0$, we have $\dot\phi_i = -K\phi_i - \beta$, whose solution is $\phi_i(t)=\big(\phi_{i,0} + \frac{\beta}{K}\big)\e^{-Kt} - \frac{\beta}{K}$. Thus, $\phi_i$ vanishes at $t_i := \frac{1}{K}\ln\big(\frac{K\phi_{i,0} + \beta}{\beta}\big)$.
By symmetry, we have that if $\phi_{i,0}<0$, the reaching time is $t_i := \frac{1}{K}\ln\big(\frac{-K\phi_{i,0}+\beta}{\beta}\big)$, and if $\phi_{i,0} = 0$, then $t_i = 0$.
Moreover, $\{\phi_i=0\}$ is forward invariant, since $\dot\phi_i \leq - \beta<0$ whenever $\phi_i>0$ and $\dot\phi_i \geq \beta>0$ whenever $\phi_i<0$. Thus, $\phi_i(t)=0$ for all $t\in[t_i,T_{\max})$, and the statement follows by setting $t_\textup{r}=\max_i t_i$. 
\end{proof}
The reaching time $t_\textup{r}$ in Lemma~\ref{lem:finite_time} depends only on $K$, $\beta$ and the initial conditions via $|\phi_{i,0}|$. %, and not on $f$, $g$, or the evolution of $x$.
Moreover, $t_\textup{r}$ is finite because $\beta>0$. For $\beta=0$, $\phi_\epsilon(x(t),\mu(t))$ converges to $\0_m$ only exponentially. The finite-time convergence of COSMIC also serve as a feasibility certificate.
\begin{corollary}[Feasibility certificate]
    % In fact, 
    In the setting of Lemma~\ref{lem:finite_time}, if the trajectory of \cosmic~starting from some $(x_0,\mu_0)$, escapes to infinity at $T < t_\textup{r}$, then~\eqref{eq:OP_ineq} admits no strictly feasible point.
\end{corollary}
To establish convergence rates we introduce the following assumption.
\begin{assumption}[Problem structure]
\label{ass:1_eq}
Given $\map{f}{\R^n}{\R}$ and $\map{g}{\R^n}{\R^m}$, assume
\begin{enumerate}[label=\textup{(\roman*)}]
\item 
\label{ass:1_f}
the function $f$ is $L_f$-smooth and the functions $g_i$ are convex with $L_g$-Lipschitz gradient, for all $i\in\{1,\dots,m\}$;
\item 
\label{ass:3_slater}
there exists $x_{\textup{s}} \in \R^n$ such that $g(x_{\textup{s}}) < \0_m$ (Slater condition).
\end{enumerate}
\end{assumption}

Lemma~\ref{lem:finite_time} holds on $[0,T_{\max})$, but it does not guarantee that $T_{\max}>t_\textup{r}$, let alone that $T_{\max}=\infty$. Since the primal dynamics is gradient-like, the solution of~\eqref{eq:SMC-Dyn} can escape to infinity in finite time only if the dual dynamics does. The next result provides a sufficient condition to preclude this.
\begin{lem}[Global existence]
\label{lem:existence}
Let $(x(t),\mu(t))$ be a Filippov solution of~\eqref{eq:SMC-Dyn} with $\phi_\epsilon(x(0),\mu(0)) \le \0_m$.  Under Assumption~\ref{ass:1_eq}, $T_{\max}=+\infty$.
\end{lem}
\begin{proof}
We provide a sketch of the proof and give the full proof in Appendix~\ref{app:existence}.
Let $\phi_i(t) := \phi_{\epsilon,i}(x(t),\mu(t))$. Since $\phi_{i,0}\le0$, by the proof of Lemma~\ref{lem:finite_time}, $\phi_i(t)\in[\phi_{i,0},0]$ for all $t$, and $\mu_i$ can be written as an explicit function of $x$ and $\phi_i(t)$, which diverges only if $x$ approaches the boundary of $\{x : g_i(x)<-\phi_i(t)\}$.
Let $x_\textup{s}$ be a Slater point and $\mu_i^+:=\max\{\mu_i,0\}$. On any bounded time interval (i) by convexity of $g_i$ and Slater's condition, $\norm{x - x_\textup{s}}$ and $\int\sum_i\mu_i^+ dt$ are bounded and (ii) the function $V(t,x) := f(x) + \sum_{i=1}^m\big(\phi_i g_i - \tfrac{\phi_i^2+\epsilon^2}{2}\log(-\phi_i(t)-g_i(x))\big)$, grows at most by a constant times $\int\sum_i\mu_i^+ dt$, and hence stays bounded. Since $V\to+\infty$ as $x$ approaches the boundary of $\{x : g_i(x) < -\phi_i(t)\}$, $x$ stays away from that boundary. Therefore, $(x,\mu)$ remains bounded on every bounded time interval, and $T_{\max}=+\infty$.
\end{proof}
Note that $\mu(0)\le\0_m$ implies $\phi_\epsilon(x(0),\mu(0)) < \0_m$, so Lemma~\ref{lem:existence} holds whenever $\mu(0) \leq \0_m$. The lemma implies that \cosmic~can be initialized at any $x(0)\in\R^n$ together with, for example $\mu(0) = \0_m$. Hence, no projection onto the feasible set is required prior to running the algorithm.
%\begin{remark}[Initialization from an infeasible point] The proposed algorithm imposes no feasibility requirement on the initial point, and therefore does not require a projection step prior to running the algorithm.\end{remark}

\bt[Convergence of \cosmic]
\label{thm:convergence}
Under Assumption~\ref{ass:1_eq}, let $f$ be $\rho$-strongly convex, let $\beta >0$, and let $(x(t), \mu(t))$ be a Filippov solution of~\eqref{eq:SMC-Dyn} with $\phi_\epsilon(x(0),\mu(0))\le\0_m$. Then, \cosmic~converges to $(x^\star_\epsilon,\mu^\star_\epsilon)$.
In particular,
$$
\norm{x(t) - x^\star_\epsilon} \leq \e^{-\rho(t - t_\textup{r})}\norm{x(t_\textup{r}) - x^\star_\epsilon}
\quad \text{ for all } t \geq t_\textup{r},
$$
with $t_\textup{r}$ as in Lemma~\ref{lem:finite_time}, and $\mu(t)\to\mu^\star_\epsilon$ as $t\to\infty$.
\et
\begin{proof}
By Lemma~\ref{lem:existence}, $T_{\max}=+\infty$. Next, by Lemma~\ref{lem:finite_time}, $\phi_\epsilon(x(t),\mu(t))=\0_m$ $\forall t\ge t_\textup{r}$, which, by Theorem~\ref{thm:eq_KKT}.\ref{th_eqKKT:item1} implies $g(x(t))<\0_m$, $\mu(t)>\0_m$ and $\mu_i(t)g_i(x(t))=-\epsilon^2/2$, $\forall i$, and $ \forall t\ge t_\textup{r}$.
Substituting $\mu_i=-\epsilon^2/(2g_i(x))$ into~\eqref{eq:SMC-Dyn} yields $\dot x = -\nabla f(x) + \frac{\epsilon^2}{2}J_g^\top \1_m\oslash g(x) = -\nabla B_\epsilon(x)$, where $B_\epsilon(x)= f(x)-\frac{\epsilon^2}{2}\sum_{i=1}^m\log{(-g_i(x))}$. That is, $x$ follows the gradient flow of the logarithmic barrier function of Fiacco–McCormick~\citep{AF-GM:1990} defined on the feasible set $\mathcal{F}$.
By assumption, $B_\epsilon(x)$ is $\rho$-strongly convex on the nonempty set $\setdef{x\in\R^n}{g(x)<\0_m}$, and its sublevel sets are compact.
Then, $B_\epsilon$ has a unique minimizer $x^\star_\epsilon$ and $\nabla B_\epsilon(x^\star_\epsilon)=\0_n$. Along the solution of~\eqref{eq:SMC-Dyn},
$
\tfrac{d}{d t}\|x-x^\star_\epsilon\|^2 
%=-2(\nabla B_\epsilon(x)-\nabla B_\epsilon(x^\star_\epsilon))^\top(x-x^\star_\epsilon)
\leq -2\rho\|x-x^\star_\epsilon\|^2,
$
hence
$\|x(t)-x^\star_\epsilon\|\leq \e^{-\rho(t-t_\textup{r})}\|x(t_\textup{r})-x^\star_\epsilon\|$.
Finally, since $g(x^\star_\epsilon)<\0_m$, continuity of $g$ gives
$\mu(t)=-\tfrac{\epsilon^2}{2}\1_m\oslash g(x(t))\to\mu^\star_\epsilon$. 
\end{proof}
Under the assumptions of Theorem~\ref{thm:convergence}, problem~\eqref{eq:OP_ineq} has a unique minimizer and satisfies Slater's condition. By Theorem~\ref{thm:eq_KKT}.\ref{th_eqKKT:item4}, from any $x(0)\in\R^n$ and, e.g., $\mu(0)=\0_m$, \cosmic~reaches the strictly feasible region in finite time and converges exponentially to $x^\star_\epsilon$ with $f(x^\star_\epsilon) - f^\star \leq m\epsilon^2/2$.

\paragraph{Relation with barrier methods.}
Logarithmic barrier methods (see, e.g.,~\citep[Ch.~11]{SB-LV:04}) approximate~\eqref{eq:OP_ineq} by the family of problems
\beq
\label{eq:barrier}
\min_{x\in\R^n} f(x) - \frac{1}{p}\sum_{i=1}^m \log\big(-g_i(x)\big), \quad p > 0,
\eeq
whose solutions $x^\star_p$ define the central path. Setting $\lambda_{p,i}^\star := - 1/(p g_i(x^\star_p))$, the optimality conditions of~\eqref{eq:barrier} can be written as the perturbed KKT system
\beq\label{eq:KKT_interior_point}
\nabla f(x)+J_g(x)^\top\lambda = \0_n,\quad g(x) < \0_m, \quad \lambda>\0_m,\quad \lambda_i g_i(x) = - \dfrac{1}{p},\quad i\in\{1,\dots,m\},
\eeq
and, in the convex case, $f(x^\star_p)-f^\star\!\leq\!m/p$~\citep[Sec.~11.2]{SB-LV:04}.
By Theorem~\ref{thm:eq_KKT}.\ref{th_eqKKT:item1}, \cosmic's equilibria are the solutions of~\eqref{eq:KKT_interior_point} with $p\!=\!2/\epsilon^2$. %i.e., the stationary points of~\eqref{eq:barrier} together with the corresponding multipliers. 
In the convex case, the $x$-components of $\Omega_\epsilon$ are the minimizers of~\eqref{eq:barrier} with $p = 2/\epsilon^2$, %i.e., points of the central path, 
and the bound in Theorem~\ref{thm:eq_KKT}.\ref{th_eqKKT:item4} recovers the classical duality gap $m/p$. This provides a control-theoretic interpretation of the central path. Specifically, its points are the equilibria of~\eqref{eq:SMC-Dyn}. However, the barrier function in~\eqref{eq:barrier} is defined on the strict interior of the feasible set, so classical barrier methods must be initialized at a strictly feasible point, which in general requires solving an auxiliary phase~\citep[Sec.~11.4]{SB-LV:04}. In \cosmic, instead, $\mu$ is a state of the dynamics driven by the SMC law~\eqref{eq:Sliding law}, and strict feasibility is attained in finite time, once the trajectories reach the sliding manifold $\mathcal S$. \label{sec: results}
\subsection{Algorithm Design}
% Discrete time algorithm + conv results?
% Trade-off
% Add references to modifications
% However, we must note the potential numerical singularities that arise in the discrete-time implementation of COSMIC. 
% Although $D_\epsilon^{-1}$ exists at the active constraints, it is increasingly ill-conditioned as the value of $\epsilon$ diminishes. To remedy this, we introduce a regularizing constant $\delta_D>0$ to ensure that $D_\epsilon$ is well conditioned. 
% We introduce nominal values for the additional hyperparameters, with $c_\phi \approx 0.85$ representing the shrinking gate based on which $\epsilon_i$ decays at a rate proportional to $r_\epsilon \approx 1.02$.

%The results of Section~\ref{sec: conv} extend verbatim to a vector $\epsilon$, with $m\epsilon^2/2$ replaced by $\norm{\epsilon}^2/2$.
  Algorithm~\ref{alg:COSMO_step} is one step of the Euler discretized equivalents of~\eqref{eq:SMC-Dyn} with two additional modifications. First, on the sliding manifold, $D_{\epsilon,i} = 1-{\mu_i}/{\sqrt{\mu_i^2+\epsilon_i^2}}$ for active constraints, so $D_\epsilon^{-1}$, though well defined, becomes ill-conditioned as $\epsilon$ decreases. We therefore replace $D_{\epsilon,i}^{-1}$ by $(D_{\epsilon,i}^2+\delta_D^2)^{-1/2}\le 1/\delta_D$, with $\delta_D>0$. Second, to reduce chattering, we use the boundary-layer approach for SMCs and replace $\operatorname{sign}(\cdot)$ by the component-wise saturation function~\citep{HK:02}.
  Theorem~\ref{thm:eq_KKT} suggests decreasing $\epsilon$ to approach a KKT point. We implement this via a vectorized $\epsilon$-homotopy that can be seen in Algorithm~\ref{alg:COSMO_solver_v2}, where  $\epsilon_i$ is divided by $r_\epsilon>1$, with a lower bound $\epsilon_{\min}>0$ set to machine precision (in our case $\epsilon_{\min} = 10^{-16}$). Specifically, $\epsilon_i$ decays when $|\phi_{\epsilon,i}|<c_\phi\,\epsilon_i$ with $c_\phi\in\,]0,1[$, as in path-following interior-point methods. Using a separate $\epsilon_i$ per constraint lets $\epsilon_i$ decay freely for inactive constraints, without being limited by the active constraints, for which a small $\epsilon_i$ causes numerical stiffness.
  Unless stated otherwise, we use $c_\phi = 0.85$ and $r_\epsilon = 1.05$, which are values that we find to work well in practice.
 Since \(D_\epsilon\) is diagonal, its inversion is computationally inexpensive and is implemented as a component-wise division of vectors.
\begin{algorithm}[h]
\caption{COSMIC \protect\Step{} function}
\label{alg:COSMO_step}
\textbf{Input:}~$x_k,\ \mu_k,\ \epsilon_k$; gains~$K,\beta>0$; step size~$dt$; saturation width~$\epsilon_{\mathrm{sat}}$; damping floor~$\delta_{D}$\\
% \textbf{Output:}~$x_{k+1},\ \mu_{k+1},\ \phi_k$\\,\ $\delta_{D} \gets c_{D}Kdt$
\ $x_{k+1} \gets x_k - dt(\nabla f(x_k) + J_g(x_k)^{\top}\mu_k)$\\
$s_k \gets \sqrt{g(x_k)^{2}+\mu_k^{2}+(\epsilon_{k})^{2}}$, \ $\phi_k \gets \mu_k- g(x_k) - s_k$, \ $d_k \gets \mathbf{1} - \mu_k \oslash s_k$\\ %d_{k,i} \in [0,1]
$D_{\epsilon,k} \gets \operatorname{diag}\big(\sqrt{d_k^{2}+\delta_{D}^{2}}\big)^{-1}$, \ $G_k \gets \operatorname{diag}\big(1+g(x_k)\oslash s_k\big)$\\
$\mu_{k+1} \gets \mu_k-dt\,D_{\epsilon,k}\Big[K\phi_k+\beta\,\mathrm{sat}\big(\phi_k,\epsilon_{\mathrm{sat}}\big)+ G_k J_g(x_k)(\nabla f(x_k) + J_g(x_k)^{\top}\mu_k)\Big]$\\
\Ret $x_{k+1},\ \mu_{k+1},\ \phi_k$
\end{algorithm}
 
\begin{algorithm}[h]
\caption{COSMIC with smooth conditioning floor and safeguarded $\epsilon$-homotopy}
\label{alg:COSMO_solver_v2}
\textbf{Input:}~$x_0$;~$\mu_0$;~$\epsilon_{\min}$; initial base smoothing~$\epsilon_0$; horizon~$N$; stationarity, feasibility, and complementarity tolerances~$(\tau_{\mathrm{s}},\tau_{\mathrm{f}},\tau_{\mathrm{c}})$; conditioning margin~$c_D$; shrink safeguard~$c_{\phi}$; homotopy rate~$r_\epsilon$\\
% \textbf{Output:}~$x_N,\ \mu_N$\\
\While{$ \big\|\nabla f(x_{k})+J_g(x_{k})^\top\mu_{k}\big\|\!>\!\tau_{\mathrm{s}}$, $\max\big(g(x_{k}),\,0\big)\!>\!\tau_{\mathrm{f}}$, $\max_i\big|\mu_{k,i}\,g_i(x_{k})\big|\!>\!\tau_{\mathrm{c}}$, \textbf{and} {$k\!<\!N$}}{
    $x_{k+1},\ \mu_{k+1},\ \phi_k \gets$ \Step{$x_k,\ \mu_k,\ \epsilon_{k}$} \hfill (Algorithm~\ref{alg:COSMO_step})\\
    $\bm{m} \gets \1_{\big[\,|\phi_{k,i}| < c_{\phi}\,\epsilon_{k,i}\,\big]_{i=1}^m}$, \ $\epsilon_{\mathrm{shrunk}} \gets \max\big(\epsilon_{k}/r_\epsilon,\ \epsilon_{\min}\big)$, \ $\epsilon_{k+1} \gets \bm{m}\odot\epsilon_{\mathrm{shrunk}} + (\mathbf{1}-\bm{m})\odot\epsilon_{k}$\\
    $r_{\mathrm{s}} \gets \big\|\nabla f(x_{k+1})+J_g(x_{k+1})^\top\mu_{k+1}\big\|$, \ $r_{\mathrm{f}} \gets \max\big(g(x_{k+1}),\,0\big)$, \ $r_{\mathrm{c}} \gets \max_i\big|\mu_{k+1,i}\,g_i(x_{k+1})\big|$
    % \If{$r_{\mathrm{s}}<\tau_{\mathrm{s}}$ \textbf{and} $r_{\mathrm{f}}<\tau_{\mathrm{f}}$ \textbf{and} $r_{\mathrm{c}}<\tau_{\mathrm{c}}$}{
    %     \Brk
    % }
    % \lIf{$r_{\mathrm{s}}<\tau_{\mathrm{s}}$ \textbf{and} $r_{\mathrm{f}}<\tau_{\mathrm{f}}$ \textbf{and} $r_{\mathrm{c}}<\tau_{\mathrm{c}}$}{\Brk}
    }
\Ret $x_N,\ \mu_N$
\end{algorithm}

% \subsection{Infeasibility Certificates}

% We characterize certificates using the upper bound on the reaching time $t_r$, which can be used to indicate the infeasibility of both convex and nonconvex optimization problems. 
 \label{sec: algorithm}
\section{Numerical Examples}
We now validate \cosmic~on SVM optimization for the Iris Setosa-vs-Versicolor dataset, and on image completion via ICNNs. 
We compare \cosmic~to Feedback linearization (FL)~\cite{RZ-AR-JS-NL:25}, Adam~\citep{DK-JB:2014} within an Augmented Lagrangian method (Adam-AL), and Adam applied to a penalty function (Adam-PEN)~\citep{lu2021physics}.
\begin{remark}
% We translate the provided implementation of the FL algorithm in~\cite{RZ-AR-JS-NL:25} from NumPy to JAX to streamline the simulations and leverage JAX’s just-in-time compilation capabilities. The algorithmic structure and the use of Gurobi to solve the QP to obtain the Lagrange multipliers, are retained, which is the primary contributor to FL’s substantial computational cost.
We re-implement the FL algorithm from~\cite{RZ-AR-JS-NL:25} in JAX with just-in-time compilation to speed up simulations. The algorithmic structure is unchanged, including the Gurobi quadratic program (QP) solver for the multipliers, which increases FL’s computational cost.
\end{remark}
% \paragraph{Rosen-Suzuki Problem.}
% We consider the nonconvex Rosen--Suzuki problem with solution $x^\star = (0,1,2,-1)$~\citep{HW-KS:80}:
% \begin{equation}\label{OP:Rosen-Suzuki}
% \begin{aligned}
% \min_{x\in\R^4} & f(x) = x_1^2 + x_2^2 + 2x_3^2 + x_4^2 - 5x_1 - 5x_2 - 21x_3 + 7x_4, \\
% \text{s.t. }&g_1(x) = x_1^2 + x_2^2 + x_3^2 + x_4^2 + x_1 - x_2 + x_3 - x_4\le 8,\\
% &g_2(x) = x_1^2 + 2x_2^2 + x_3^2 + 2x_4^2 - x_1 - x_4 \le 10,\\
% &g_3(x) = 2x_1^2 + x_2^2 + x_3^2 + 2x_1 - x_2 - x_4 \le 5.
% \end{aligned}
% \end{equation} 
% We sample $N = 10$ random initial points drawn uniformly over $[-5,5]^4$, and present the results for COSMIC, FBL-M \citep{RZ-AR-JS-NL:25}, and KKT-SMC \citep{AZ-KY:2026} in Figure~\ref{fig: Rosen-Suzuki}.
% \begin{figure}
%     \centering
%     \includegraphics[width=0.5\linewidth]{}
%     \caption{Optimization of the Rosen-Suzuki Problem}
%     \label{fig: Rosen-Suzuki}
% \end{figure}
% Though KKT-SMC is able to attain the minimum 
\begin{wrapfigure}{H}{0.6\linewidth}
    \centering
    \includegraphics[width =0.95\linewidth, page = 2]{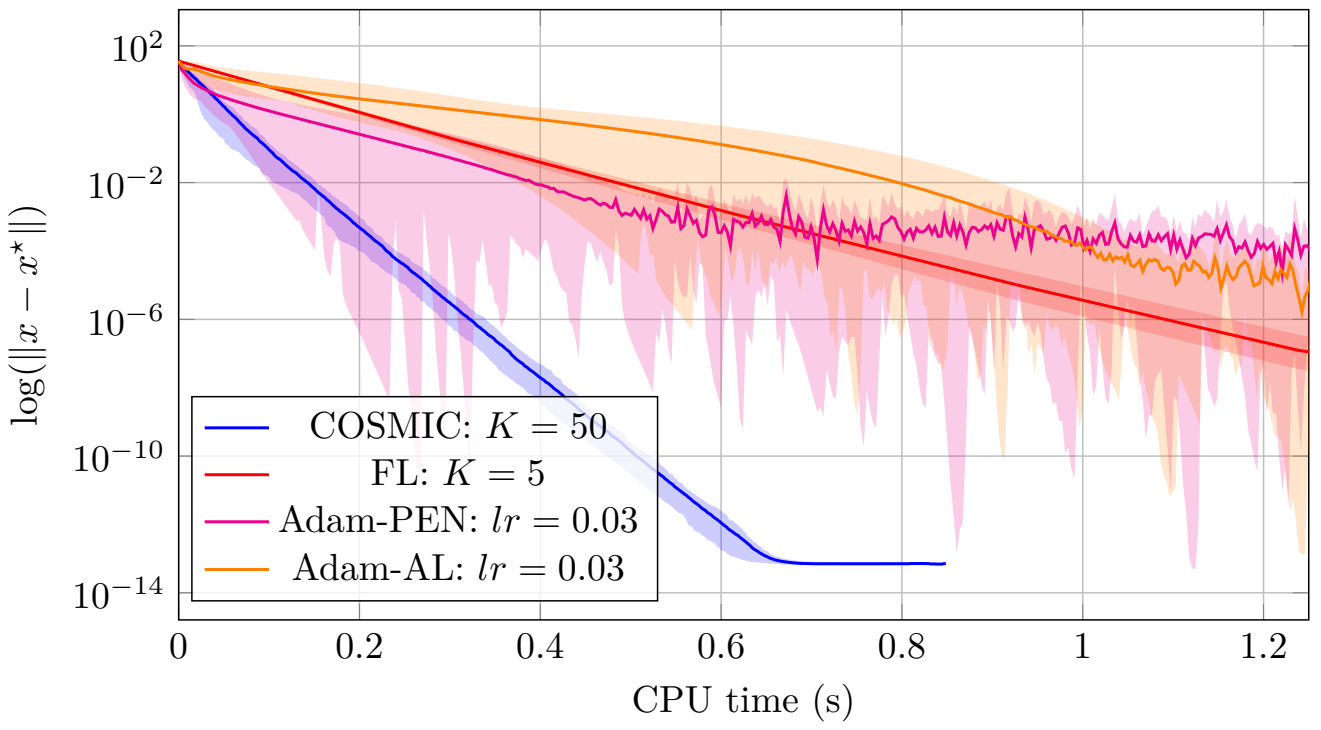}
    \caption{Optimization of soft-margin SVM for Iris Dataset. The solid line is the mean error over the $10$ runs, and the shaded band corresponds to their minimum and maximum.}
    \label{fig: SVM-SM}
\end{wrapfigure}
\paragraph{Iris Dataset SVM.} We consider a support vector machine (SVM) benchmark on real data~\citep{iris-dataset}. Training a soft-margin SVM is equivalent to $\min_{w,b,\xi} \frac{1}{2}\|w\|^2 + (1/2)\|\xi\|^2 ~\text{s.t. } 1 - y_i(w^\top a_i + b) - \xi_i \leq 0$, $i=1,\dots,N_\text{s},
$ where $a_i\in\mathbb{R}^{d}$ are the samples with labels $y_i\in\{-1,+1\}$, and $N_\text{s}$ is the number of samples. We use the full dataset ($N_\text{s} = 150$) with $d=4$, and train a soft-margin SVM to separate the Virginica species from the other two. The resulting formulation has $n = 155$ optimization variables, and the optimal solution for both problems is computed using SLSQP and then refined by Newton's method on the active-set KKT system. Each method is run from $10$ different initial conditions, with components sampled independently and uniformly from $[-1,1]$. Since, with the tested parameters, the second-order KKT-SMC method fails to converge from any initial condition, we defer its comparison with \cosmic~to Appendix~\ref{app:SVM}, on a hard-margin SVM problem with a smaller dataset. Figure~\ref{fig: SVM-SM} plots the $\log$ of the error of FL, Adam-AL, Adam-PEN, and \cosmic~with respect to CPU time. \cosmic~significantly outperforms the other methods, reaching the highest accuracy in the shortest CPU time. Although FL reaches a final accuracy of order $10^{-9}$, its per-iteration cost is high, averaging $6.14$\,ms versus $0.042$\,ms for \cosmic, $0.059$\,ms for Adam-PEN, and $0.066$\,ms for Adam-AL.

\begin{figure}[!h]
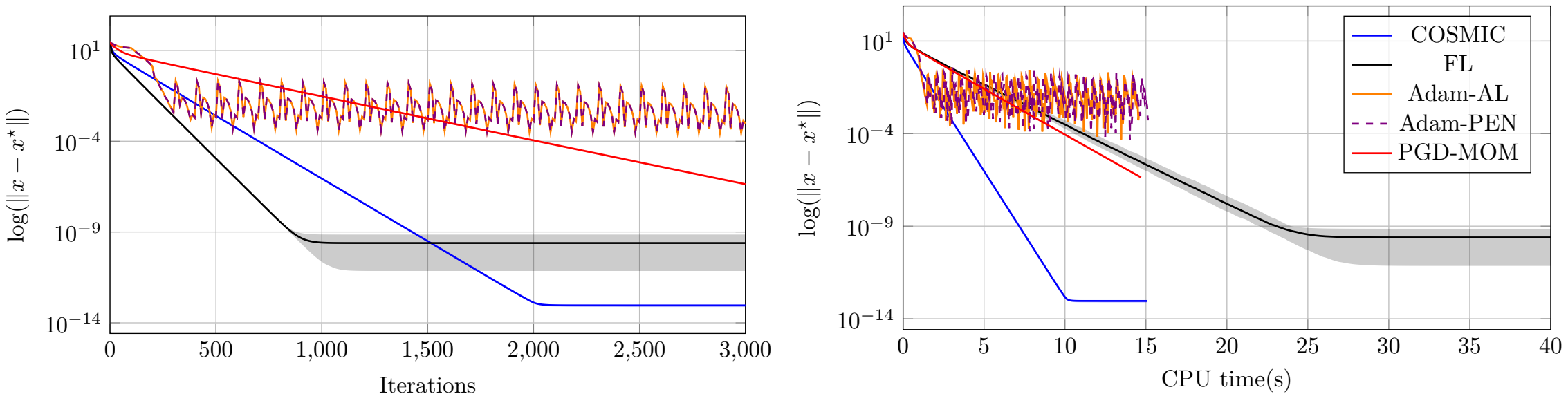

  \centering
  \includegraphics[width =0.49\linewidth, page = 3]{Figures/L4DC-COSMIC-Figures.pdf}
  \includegraphics[width =0.49\linewidth, page = 4]{Figures/L4DC-COSMIC-Figures.pdf}
  \caption{Error convergence with respect to (a) iterations and (b) CPU time. The solid line shows the mean error over the $3$ runs, and the shaded band spans their minimum and maximum. For visualization purposes, the CPU-time range is limited to $40 $s.}
  \label{fig:icnn}
\end{figure}

\paragraph{Image Completion using ICNN} We consider a higher-dimensional problem: completing the right half of a $64\times64$ face image from its left half with a learned convex energy given by a partially input convex neural network (PICNN)~\citep{BA-LC-JK:2017}. The completed half ($n=2048$ pixels) is obtained by minimizing the PICNN output subject to $m=4096$ linear inequality constraints that enforce the valid pixel range $[0,1]$. 
Figure~\ref{fig:icnn} compares \cosmic~with FL, Adam-AL and Adam-PEN, and with projected gradient descent with momentum (PGD-MOM), the inference method of~\citep{BA-LC-JK:2017}. 
All methods are run for $3000$ iterations on the same $3$ test images and $2$ initial infeasible points per image drawn uniformly from $[-0.5,1.5]^n$. 
% with the exception of FL which is run for $1000$ iterations since it solves an quadratic program at every iteration. 
The CPU time is measured per iteration on the compiled solvers, excluding metric logging; for FL it includes the QP solve in Gurobi. The hyperparameters for PGD-MOM are adopted from~\cite{BA-LC-JK:2017}, while those of the other methods are tuned separately via a parameter sweep for best performance, and listed in~Table~\ref{tab:hp_icnn} in Appendix~\ref{app:icnn}.
Figure~\ref{fig:icnn} shows that \cosmic~reaches the highest accuracy in the shortest CPU time, and has better and comparable performance to adaptive momentum based methods. Although While FL achieves high accuracy, the QP solved at each iteration results in an average iteration time of approximately $24$\,ms, considerably higher than the $4.8$--$5.5$\,ms of the other methods, with a total runtime exceeding $80$\,s for $3000$ iterations. Adam-AL, Adam-PEN, and PGD-MOM converge more slowly and do not reach the accuracy of \cosmic~within the same budget.
 \label{sec: numerical_results}
\section{Conclusion}
% In this paper, we introduce continuous-time gradient-flow dynamics derived by designing the Lagrange multiplier dynamics of inequality-constrained optimization problems using SMC. We refer to the dynamics as \cosmic~for Constrained Optimization via Sliding Mode for Inequality Constraints. We provide a control-theoretic perspective on log barrier methods, whcih are equivalent to COSMIC after its finite time convergence to the smoothed FB function and address the numerical stiffness arising from the implementation to obtain an SMC-based solver for inequality-constrained optimization. We validate COSMIC across various machine learning applications, demonstrating its ability to solve optimization problems accurately while scaling effectively to larger problem instances, which is necessary for practical machine learning applications. Furthermore, COSMIC's finite-time convergence property provides a means of certifying infeasibility which we will explore in more details in future work. 
% COSMIC not only achieves finite-time convergence but also demonstrates improved performance compared with existing optimization approaches, including ADAM with an augmented Lagrangian formulation.
We introduce \cosmic~\eqref{eq:SMC-Dyn}, a continuous-time solver for inequality-constrained optimization in which the Lagrange multiplier dynamics are designed via SMC. Its trajectories reach, in finite time, the sliding manifold where the smoothed FB complementarity condition are satisfied. On this manifold, \cosmic~is the gradient flow of a log-barrier function, which yields a control-theoretic interpretation of barrier methods.
We show that the equilibria of~\eqref{eq:SMC-Dyn} are exactly the points satisfying stationarity and the smoothed FB condition, and that, as the smoothing parameter $\epsilon$ decays, every limit point of a sequence of equilibria is a KKT point of the original problem. For convex problems satisfying Slater's condition, the suboptimality is bounded by $m\epsilon^2/2$. For strongly convex costs, from any initial primal point, the trajectories become strictly feasible in finite time and then converge exponentially to the minimizer of the associated barrier problem. Experiments on SVM classification and image completion via an ICNN mapping confirm that \cosmic~solves these problems accurately and scales up without significant increase in computational time per-iteration. 
% Future directions include the exploration of \cosmic's finite-time convergence to the sliding manifold to certifying infeasibility. 
Future directions include the extension of the convergence analysis beyond strong convexity, analysis of discrete-time implementations, and extensions to accommodate equality constraints. \label{sec: conclusion}
\subsection*{Generative AI Use}
Generative AI tools were used for code generation, editing, and proof checking. The authors reviewed and verified all AI-assisted content and take responsibility for the final content of the paper.
\acks{This work was supported as a part of NCCR Automation, a National Centre of Competence in Research, funded by the Swiss National Science Foundation (grant number 51NF40\_225155).
}
\bibliography{main}

@article{AF:92,
  title = {A special {N}ewton-type optimization method},
  volume = {24},
  ISSN = {1029-4945},
  DOI = {10.1080/02331939208843795},
  number = {3-4},
  journal = {Optimization},
  author = {A. Fischer},
  year = {1992},
  month = Jan,
  pages = {269–284}
}

@article{CK:96,
  title = {Some Noninterior Continuation Methods for Linear Complementarity Problems},
  volume = {17},
  ISSN = {1095-7162},
  DOI = {10.1137/s0895479894273134},
  number = {4},
  journal = {SIAM Journal on Matrix Analysis and Applications},
  author = {C. Kanzow},
  year = {1996},
  month = Oct,
  pages = {851–868}
}

@Book{SB-LV:04,
  author = {S. Boyd and L. Vandenberghe},
  title = {Convex Optimization},
  publisher = {Cambridge University Press},
  year = 2004,
  isbn = 0521833787,
}

@book{HK:02,
  title={Nonlinear Systems},
  author={H. K. Khalil},
  isbn={9780130673893},
  lccn={95045804},
  series={Pearson Education},
  year={2002},
  publisher={Prentice Hall}
}

@article{BD-SP-SK-TD:2026,
  title={Finite and fixed-time feedback-based continuous-time optimization},
  author={B. Diana and S. Pandey and S. Kamal and T. N. Dinh},
  journal={Automatica},
  volume={183},
  pages={112569},
  year={2026},
  DOI = {10.1016/j.automatica.2025.112569}
}

@article{FD-ZH-GB-SB-JL-MM:24,
  title = {Toward a Systems Theory of Algorithms},
  author = {F. D\"{o}rfler and Z. He and G. Belgioioso and S. Bolognani and J. Lygeros and M. Muehlebach},
  volume = {8},
  ISSN = {2475-1456},
  DOI = {10.1109/lcsys.2024.3406943},
  journal = {IEEE Control Systems Letters},
  year = {2024},
  pages = {1198--1210}
}

@inproceedings{RZ-AR-JS-NL:25,
  title = {Constrained Optimization From a Control Perspective via Feedback Linearization},
  author={R. Zhang and A. Raghunathan and J. Shamma and N. Li},
  DOI = {10.52202/085713-1873},
  booktitle = {Advances in Neural Information Processing Systems 38},
  year = {2025},
  pages = {62328--62364}
}

@misc{AZ-KY:2026,
	title = {{KKT} {Residual} {Sliding} {Manifolds}: {A} {Projection}-{Free} {Structured} {Finite}-{Time} {Solver} for {Nonlinear} {Constrained} {Optimisation}},
	author = {A. Zenati and K. Youcef-Toumi},
    journal = {arXiv preprint arXiv:2610.15545},
	year = {2026}
}

@book{AF-GM:1990,
  title={Nonlinear programming: sequential unconstrained minimization techniques},
  author={A. V. Fiacco and G. P. McCormick},
  year={1990},
  publisher={SIAM}
}

@article{VC-RH-EB-JL:2026,
  title={{SPPID}: Saddle-Point {PID} for Constrained Optimization},
  author={V. Centorrino and R. Hoteit and E. C. Balta and J. Lygeros},
  journal={arXiv preprint arXiv:2609.31086},
  year={2026}
}

@article{VC-SMF-SP-DR:25,
  title = {A new framework for constrained optimization via feedback control of {L}agrange multipliers},
  author={V. Cerone and S. M. Fosson and S. Pirrera and D. Regruto},
  ISSN = {2334-3303},
  DOI = {10.1109/tac.2025.3568651},
  journal = {IEEE Transactions on Automatic Control},
  publisher = {Institute of Electrical and Electronics Engineers (IEEE)},
  year = {2025},
  pages = {1--16}
}

@article{AH-ZH-SB-GH-FD:24,
  title = {Optimization algorithms as robust feedback controllers},
  volume =	 57,
  doi = {10.1016/j.arcontrol.2024.100941},
  journal =	 {Annual Reviews in Control},
  author =	 {A. Hauswirth and Z. He and S. Bolognani and G. Hug and F. D\"orfler},
  year = 2024,
  pages =	100941
}

@book{KJA-LH-HU:58,
	editor = {K. J. Arrow and L. Hurwicz and H. Uzawa},
	publisher = {Stanford University Press},
	title = {Studies in Linear and Nonlinear Programming},
	year = {1958}
}

@article{AA-JC:2023,
  title={Control-barrier-function-based design of gradient flows for constrained nonlinear programming},
  author={A. Allibhoy and J. Cort{\'e}s},
  journal={IEEE Transactions on Automatic Control},
  volume={69},
  number={6},
  pages={3499--3514},
  year={2023},
  doi = {10.1109/TAC.2023.3306492}
}

@article{GB-JC-JP-ED:2021,
  title={Time-varying optimization of LTI systems via projected primal-dual gradient flows},
  author={G. Bianchin and J. Cort{\'e}s and J. I. Poveda and E. Dall’Anese},
  journal={IEEE Transactions on Control of Network Systems},
  volume={9},
  number={1},
  pages={474--486},
  year={2021},
  DOI = {10.1109/tcns.2021.3112762}
}

@article{JC:2008,
  title={Discontinuous dynamical systems},
  author={J. Cortes},
  journal={IEEE Control Systems Magazine},
  volume={28},
  number={3},
  pages={36--73},
  year={2008},
  DOI = {10.1109/mcs.2008.919306}
}

@article{FC-WR:2020,
  title={Sign projected gradient flow: A continuous-time approach to convex optimization with linear equality constraints},
  author={F. Chen and W. Ren},
  journal={Automatica},
  volume={120},
  pages={109156},
  year={2020},
  DOI = {10.1016/j.automatica.2020.109156}
}

@inproceedings{BA-LC-JK:2017,
  title={Input convex neural networks},
  author={B. Amos and L. Xu and J. Z. Kolter},
  booktitle={International Conference on Machine Learning},
  pages={146--155},
  year={2017},
  organization={PMLR}
}

@article{DK-JB:2014,
  author = {D. P. Kingma and J. Ba},
  title = {{Adam}: {A} Method for Stochastic Optimization},
  journal = {arXiv preprint ARXIV:1412.6980},
  year = {2014}
}

@article{CDC:VC-RH-EB-JL:2026,
  title={A Unified Control-Theoretic Framework for Saddle-Point Dynamics in Constrained Optimization},
  author={V. Centorrino and R. Hoteit and E.C. Balta and J. Lygeros},
  journal = {arXiv preprint arXiv:2604.09252},
  year={2026}
}

@misc{RTR:70,
  title={Convex Analysis},
  author={R. T. Rockafellar},
  year={1970},
  publisher={Princeton University Press}
}

@misc{iris-dataset,
  doi = {10.24432/C56C76},
  author = {R. A. Fisher},
  title = {Iris Dataset},
  publisher = {UCI Machine Learning Repository},
  year = {1936}
}

@article{lu2021physics,
  title={Physics-informed neural networks with hard constraints for inverse design},
  author={L. Lu and R. Pestourie and W. Yao and Z. Wang and F. Verdugo and S. G. Johnson},
  journal={SIAM Journal on Scientific Computing},
  volume={43},
  number={6},
  pages={B1105--B1132},
  year={2021},
  publisher={SIAM}
}
%this appendix is not in the template
\appendix
% This forces the prefix to be the Alphabetic version of the section counter
\renewcommand{\thelem}{\Alph{section}.\arabic{lem}}
\newpage
\section{Extended proof}
\label{app:existence}
\begin{proof}[Proof of Lemma~\ref{lem:existence}]
The right-hand side of~\eqref{eq:SMC-Dyn} is measurable and locally essentially bounded. Hence a Filippov solution exists from every initial condition~\citep[Proposition~3]{JC:2008}, and, a maximal Filippov solution does not exist in finite time if and only if it leaves every compact set.
Let $(x(t),\mu(t))$ be a maximal Filippov solution with $\phi_\epsilon(x(0),\mu(0))\le\0_m$, defined on $[0,T_{\max})$, and assume, for the sake of contradiction, that $T_{\max}<\infty$. In what follows, $c>0$ denotes a generic constant which may depend on $T_{\max}$, the initial condition, and the problem data. The value of $c$ may change from one occurrence to the next, even within the same line, but we keep the same symbol for simplicity of notation.

Since $\phi_{i,0}\le 0$, the proof of Lemma~\ref{lem:finite_time} gives $\phi_i(t) := \phi_{\epsilon,i}(x(t),\mu(t)) \in [\phi_{i,0},0]$ and $0\le\dot\phi_i\le K|\phi_{i,0}|+\beta$ for almost every $t\in[0,T_{\max})$.
Moreover, $d_i := -\phi_i - g_i(x) = \sqrt{g_i^2(x) + \mu_i^2 + \epsilon^2} - \mu_i$, where $d_i>0$. Expanding $(d_i+\mu_i)^2 = g_i^2(x) + \mu_i^2 + \epsilon^2$ and substituting $g_i^2(x) = \phi_i^2 + 2\phi_i d_i + d_i^2$ yields
\beq
\label{eq:app_mu}
\mu_i = \phi_i + \frac{\phi_i^2+\epsilon^2}{2d_i},
\quad \text{equivalently,} \quad
(\mu_i-\phi_i)\,d_i = \frac{\phi_i^2+\epsilon^2}{2}.
\eeq
Therefore, $\mu_i > \phi_i \geq \phi_{i,0}$, and $\mu_i$ can diverge only if $d_i \to 0$ for some $i$. It thus suffices to show that $x$ is bounded and that each $d_i$ is bounded away from zero.
Next, let $x_\textup{s}$ be a Slater point and $\gamma := \min_{j}\big(-g_j(x_\textup{s})\big)>0$. Then, by~\eqref{eq:SMC-Dyn}
$$
\dfrac{d}{dt}\dfrac{1}{2}\norm{x-x_\textup{s}}^2 = - \nabla f(x)^\top (x - x_\textup{s}) - \sum_{i=1}^m\phi_i\nabla g_i(x)^\top(x-x_\textup{s}) - \sum_{i=1}^m(\mu_i - \phi_i)\nabla g_i(x)^\top(x-x_\textup{s}).
$$
By the Cauchy--Schwarz inequality and since $\nabla f$ is $L_f$-Lipschitz, it holds $-\nabla f(x)^\top(x - x_\textup{s}) \leq \|\nabla f(x)\| \|x-x_\textup{s}\| \leq \big(\|\nabla f(x_\textup{s})\| + L_f\|x-x_\textup{s}\|\big)\|x-x_\textup{s}\| \le c (1+\|x-x_\textup{s}\|^2)$. Moreover, since $|\phi_i|\le|\phi_{i,0}|$ and $\nabla g_i$ is $L_g$-Lipschitz, the second sum is bounded by $c\,(1+\|x-x_\textup{s}\|^2)$. Since $\mu_i-\phi_i>0$ and $g_i(x)=-\phi_i-d_i\ge-d_i$, convexity of $g_i$ and~\eqref{eq:app_mu} yield $(\mu_i-\phi_i)\nabla g_i(x)^\top(x-x_\textup{s})\ge(\mu_i-\phi_i)\big(g_i(x)-g_i(x_\textup{s})\big) \geq \gamma(\mu_i-\phi_i)-c$. Thus, 
$$
\dfrac{d}{dt}\dfrac{1}{2}\|x-x_\textup{s}\|^2 + \gamma \sum_{i=1}^m(\mu_i-\phi_i)\le c \big(1+\|x-x_\textup{s}\|^2\big),
$$
 where $c$ subsumes the scaling constants. Moreover, since $T_{\max}<\infty$, Gr\"onwall's inequality gives $\|x(t)-x_\textup{s}\|\le c$ and $\int_0^t\sum_i(\mu_i-\phi_i) d\tau \leq c$ on $[0,T_{\max})$.

By~\eqref{eq:app_mu}, $\dot x=-\nabla_x V(t,x)$, where
$$
V(t,x) := f(x)+\sum_{i=1}^m\Big(\phi_i(t)\,g_i(x)-\frac{\phi_i(t)^2+\epsilon^2}{2}\log\big(-\phi_i(t)-g_i(x)\big)\Big).
$$
Along the solution, and by using $\dot x = -\nabla_x V(t,x)$, $d_i = -\phi_i - g_i(x)$, and~\eqref{eq:app_mu}, we have
\begin{align*}
\frac{d}{dt}V(t,x(t))
&= \nabla_x V(t,x)^\top\dot x + \sum_{i=1}^m \frac{\partial V}{\partial\phi_i}(t,x)\,\dot\phi_i\\
&= -\norm{\nabla_x V(t,x)}^2 + \sum_{i=1}^m \dot\phi_i\Big(g_i(x) - \phi_i\log d_i + \frac{\phi_i^2+\epsilon^2}{2\,d_i}\Big)\\
&\le \sum_{i=1}^m \dot\phi_i\big(g_i(x) - \phi_i\log d_i + \mu_i-\phi_i\big)\\
&\le \sum_{i=1}^m \dot\phi_i\big(c + \mu_i-\phi_i\big)
\le c\Big(1+\sum_{i=1}^m(\mu_i-\phi_i)\Big),
\end{align*}
where we used that $x$ is bounded, so that $|g_i(x)|\le c$ and $d_i\le|\phi_{i,0}|+|g_i(x)|\le c$, that $-\phi_i\log d_i\le|\phi_{i,0}|\max\{0,\log c\}\le c$ since $-\phi_i\ge0$, and that $\mu_i-\phi_i>0$ and $0\le\dot\phi_i\le c$.
Integrating gives $V(t,x(t))\le  c$. On the other hand, since $x$ is bounded and $d_j\le c$ for all $j$, all terms of $V$ except the $i$-th logarithm are bounded below, and the latter is at least $\frac{\epsilon^2}{2}|\log d_i|$ whenever $d_i\le1$. Hence, $d_i(t)\geq \e^{- c}$ on $[0,T_{\max})$.

By~\eqref{eq:app_mu}, $\phi_{i,0}\le\mu_i(t)\le  c$, so $(x,\mu)$ remains in a compact set on $[0,T_{\max})$, contradicting $T_{\max}<\infty$. Hence, $T_{\max} = +\infty$.
\end{proof}

\section{Numerical Experiments}
\label{app:exp}
We provide implementation details and additional results for the experiments in Section~\ref{sec: numerical_results}. Besides FL and KKT-SMC, we compare \cosmic~with two standard first-order approaches to constrained learning based on Adam, namely
%We consider alternative first-order methods to solve the constrained machine learning examples, including 
Adam with an augmented Lagrangian formulation, and Adam applied to a quadratic penalty formulation. For completeness, we review these methods.

\paragraph{Adam with Penalty Method (Adam-PEN).}
 Adam is run on the quadratic penalty function 
 $$
 F_\sigma(x) = f(x) + \dfrac{\sigma}{2}  \|\max(0, g(x))\|^2,
 $$
 where $\sigma>0$ is an increasing penalty parameter capped at $\sigma_{\max}$. 
 % Under standard regularity assumptions, the minimizer $x_\sigma$ of $F_\sigma$ satisfies $\norm{x_\sigma - x^\star} = \mathcal{O}(1/\sigma)$. %in general, however, $x_\sigma$ is infeasible, and the accuracy is limited by $\sigma$, whose increase worsens the conditioning of $F_\sigma$. 

%This leads to an error $\mathcal{O}(1/\sigma)$.
 \paragraph{Adam with an augmented Lagrangian (Adam-AL).}
As a first-order baseline, we apply the method of multipliers to the Powell--Hestenes--Rockafellar augmented Lagrangian
\begin{equation}
    \mathcal{L}_\sigma(x,\lambda)=f(x)+\frac{1}{2\sigma}\sum_{i=1}^{m}\Big(\max\{0,\lambda_i+\sigma\,g_i(x)\}^2-\lambda_i^2\Big).
\end{equation}
Each outer iteration runs $T$ Adam steps on $x$ (with $\lambda$ and $\sigma$ fixed), then updates
\begin{equation}
    \lambda\leftarrow\max\{0,\lambda+\sigma\,g(x)\},\qquad
    \sigma\leftarrow\min\{\gamma_\sigma\,\sigma,\sigma_\text{max}\},\qquad
    \alpha\leftarrow\max\{\gamma_\alpha\,\alpha,\alpha_\text{min}\},
\end{equation}
where $\alpha$ is the Adam step size, $\gamma_\sigma \geq 1$, $\gamma_\alpha \in ]0,1]$, and the Adam moments are reset after every multiplier update. The effective multiplier $\max\{0,\lambda+\sigma g(x)\}$ is used when evaluating the KKT residual.
\subsection{Hard-Margin SVM Optimization}
\label{app:SVM}
\begin{figure}[H]
    \centering
    \includegraphics[width =0.9\linewidth, page = 1]{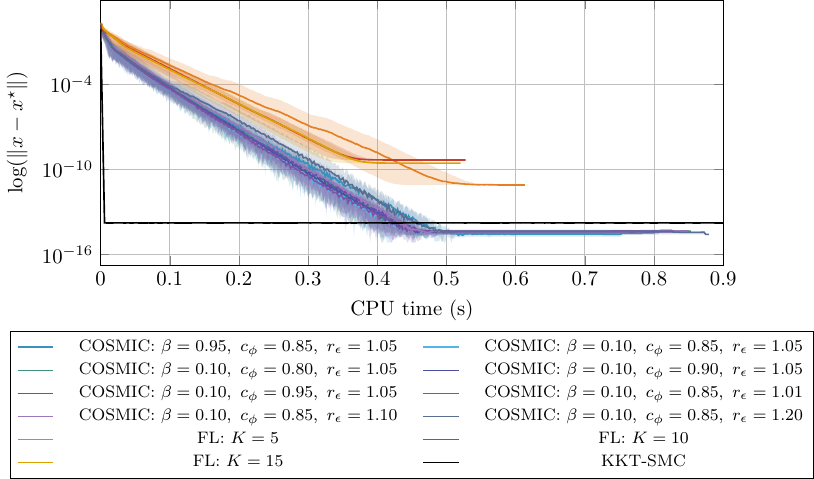}
    \caption{Optimization of Hard-Margin SVM for Iris Dataset}
    \label{fig: SVM-HM}
\end{figure}

To compare \cosmic~with KKT-SMC, we consider the hard-margin SVM problem which has less variables. In this case, training an SVM is equivalent to solving
\begin{align*}
\min_{w,b}~&\tfrac{1}{2}\|w\|_2^2\\
\text{s.t.}~&1-y_i\,(w^\top a_i+b)\le 0, \quad i=1,\dots,N_\text{s},
\end{align*}
where $a_i\in\mathbb{R}^{d}$ are the samples, $y_i\in\{-1,+1\}$ their labels, and $N_\text{s}$ is the number of samples. We use the linearly separable Iris Setosa-vs-Versicolor  dataset with $d=4$ features, which gives $n=d+1=5$ optimization variables and $N_\text{s}=100$ linear inequality constraints corresponding to only two species of the Iris flower: Setosa and Veriscolor. We aim to show the performance of \cosmic~on this lower dimensional problem and also consider its sensitivity to $c_\phi$ and $r_\epsilon$, introduced in Algorithm~\ref{alg:COSMO_solver_v2}.
All methods are run from the same $10$ initial conditions, whose components are sampled uniformly at random from $[-1,1]$, with a fixed budget of $N=20000$ iterations.
Figure~\ref{fig: SVM-HM} shows the error with respect to the optimal solution as a function of the CPU time for \cosmic, FL, and KKT-SMC, with the best-performing set of hyperparameters. The tested hyperparameters are listed in Table~\ref{tab:hp_svm}. For \cosmic, we test several values of $K$, which give comparable performance. %The parameters $c_\epsilon$ and $r_\epsilon$ are introduced in Algorithm~\ref{alg:COSMO_step}. 
For the problem formulation considered here, the performance of \cosmic~is relatively insensitive to the specific values of $c_\phi$ and $r_\epsilon$. This is illustrated in Figure~\ref{fig: SVM-HM}, where we vary $c_\phi$ over $[0.8,0.95]$ and $r_\epsilon$ over $[1.01,1.2]$.
Since the problem is low-dimensional, the second-order KKT-SMC reaches the neighborhood of the minimum within a few iterations. However, its accuracy is limited by the non-vanishing smoothing parameter of the FB function. \cosmic~requires more iterations, but attains a higher accuracy than KKT-SMC and has a lower average cost per iteration ($0.045$\,ms versus $0.105$\,ms for KKT-SMC). FL is slower and less accurate than both sliding-mode-based methods.
%We provide below the hyperparameters set for each solver for the SVM optimization problem.
\begin{table}[H]
\centering
\caption{Tested hyperparameters for SVM optimization.}
\label{tab:hp_svm}
\begin{tabular}{ll}
\hline
Method & Hyperparameters\\
\hline
\cosmic  & $d t=0.0035$, $K \in \{50,60\}$, $\beta=0.1$, $\epsilon_\text{sat}=0.01$, $\epsilon_0 = 10^{-10}$\\
FL      & $\eta=0.1$, $K \in \{5,10,15\}$\\
KKT-SMC & $\eta=0.8$, $q_x=q_c=1$, $\varepsilon_x=\varepsilon_c=0.01$, $\phi=0.1$, $\sigma_c=10^{-15}$\\
\hline
\end{tabular}
\end{table}
\subsection{Image Completion using Input Convex Neural Networks}
\paragraph{Architecture.}
Let $x\in\mathbb{R}^{64\times32}$ denote the observed left half of a face and $y\in\mathbb{R}^{n}$ the right half to be generated. The energy $F_\theta(x,\cdot):\mathbb{R}^n\to\mathbb{R}$ is a convolutional PICNN with an input path $u_i$ acting on $x$ and an output path $z_i$ acting on $y$,
\begin{equation}
\begin{aligned}
    &u_0 = x, \qquad u_{i+1} = \text{ReLU}\!\big(\tilde W_i\,u_i + \tilde b_i\big), \qquad i=0,\dots,3,\\
    &z_1 = \text{sp}\!\Big(W^y_0\big(y\circ(W^{yu}_0 x + b^y_0)\big) + W^u_0\,x + b_0\Big),\\
    &z_{i+1} = \text{sp}\!\Big(W^z_i\big(z_i\circ[W^{zu}_i u_i + b^z_i]_+\big) + W^y_i\,y + W^u_i\,u_i + b_i\Big), \qquad i=1,2,3,\\
    &F_\theta(x,y) = W^z_4\big(z_4\circ[W^{zu}_4 u_4 + b^z_4]_+\big) + W^y_4\,y + W^u_4\,u_4 + b_4 + \tfrac{\rho}{2}\|y\|_2^2 ,
\end{aligned}
\label{eq:icnn}
\end{equation}
where $\text{sp}(s)=\log(1+\e^{s})$ is the softplus function, applied element-wise, $[\cdot]_+$ is the positive part, $\circ$ is the element-wise product, and $\theta$ collects all weights and biases. We choose the network as in~\citep{BA-LC-JK:2017}.
Because $\text{sp}$ is convex and non-decreasing, the weights $W^z_i$ acting on the hidden states are element-wise non-negative, and the gates $[W^{zu}_i u_i + b^z_i]_+$ are non-negative and depend on $x$ only, every $z_i$ is convex in $y$. Hence, $F_\theta(x,y)$ is convex and smooth in $y$ for every $x$, although not necessarily convex in $x$, and the quadratic term makes it $\rho$-strongly convex with $\rho=0.05$, so~\eqref{eq:generation} has a unique minimizer.
%The map $F_\theta(x,\cdot)$ is therefore convex in $y$ for every $x$, but is not convex in $x$. The quadratic term makes it $\rho$-strongly convex with $\rho=0.05$, and the constrained problem below has a unique minimiser.

\paragraph{Training.}
The parameters $\theta$ are trained with a denoising score matching loss, such that the true right half $y$ is the minimizer of $F_\theta(x,\cdot)$: for a noisy right half $\tilde y=y+\sigma\varepsilon$, $\varepsilon\sim\mathcal{N}(0,I)$, $\sigma\sim\mathcal{U}(0.02,0.25)$,
\begin{equation}
    \min_{\theta}\; \mathbb{E}_{(x,y)\sim\mathcal{D},\,\tilde y}\!\left[\big\|\nabla_y F_\theta(x,\tilde y)-c\,(\tilde y-y)\big\|_2^2\right]
    \quad\text{s.t.}\quad W^z_i\ge0,\ \ i=1,\dots,4,
    \label{eq:training}
\end{equation}
with $c=0.25$, i.e., $F_\theta(x,\cdot)$ is fitted to the convex bowl $\tfrac{c}{2}\|\cdot-y\|_2^2$ around the true completion. Problem \eqref{eq:training} is solved with Adam (batch size $16$, learning rate $10^{-3}$ with cosine decay to $10^{-4}$, $2000$ steps), where the weights $W^z_i$ are projected onto the non-negative orthant after every step, as in~\citep{BA-LC-JK:2017}. We use the Olivetti faces ($400$ images of $40$ subjects, $64\times64$ pixels in $[0,1]$) with $36$ subjects for training, and $4$ for testing, and random horizontal flips of the full image as data augmentation.

\begin{figure}[!h]
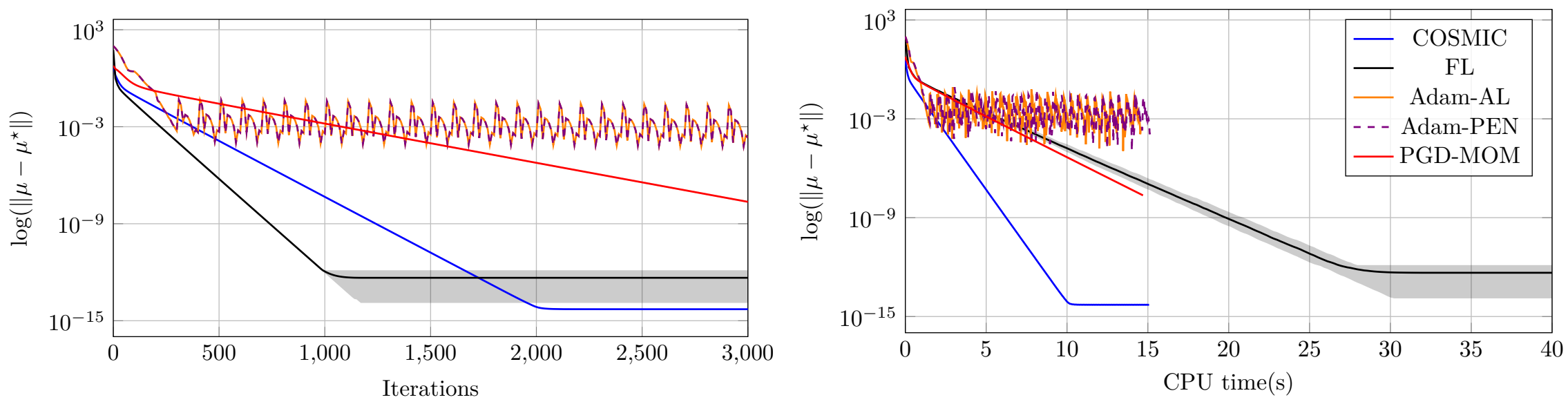

  \centering
  \includegraphics[width =0.49\linewidth, page = 5]{Figures/L4DC-COSMIC-Figures.pdf}
  \includegraphics[width =0.49\linewidth, page = 6]{Figures/L4DC-COSMIC-Figures.pdf}
  \caption{Multiplier error convergence with respect to (a) iterations and (b) CPU time. The solid line shows the mean error over the $3$ runs, and the shaded band spans their minimum and maximum. For visualization purposes, the CPU-time range is intentionally limited to $40 $s.}
  \label{fig:icnn_mult}
\end{figure}

\paragraph{Problem.}
Given the left half $x$ of a test face, the completed right half is the minimiser of the learned energy subject to the pixel range,
\begin{equation}
\begin{aligned}
     \min_{y\in\mathbb{R}^n} \ & F_\theta(x,y)\\
    \text{s.t.} \
        & \0_n\leq y\leq \1_n.
\end{aligned}
\label{eq:generation}
\end{equation}
Writing the constraints as $g(y)=Ay-b\le 0$ with
\begin{equation}
    A=\begin{bmatrix}I\\-I\end{bmatrix}\in\mathbb{R}^{m\times n},\qquad
    b=\begin{bmatrix}\1_n\\ \0_n\end{bmatrix},\qquad m=2n=4096,
\end{equation}
problem \eqref{eq:generation} becomes the strongly convex program $\min_y F_\theta(x,y)$ s.t.\ $Ay\le b$. Since the constraints are affine, the Jacobian $J_g(y)=A$ is constant and the Hessian of the Lagrangian $L(y,\mu)=F_\theta(x,y)+\mu^{\!\top}(Ay-b)$ equals $\nabla_y^2F_\theta(x,y)$. The unique primal--dual solution $(y^\star,\mu^\star)$ is characterized by the KKT conditions
\begin{equation}
    \nabla_y F_\theta(x,y^\star)+A^{\!\top}\mu^\star=\0_n,\quad
    Ay^\star-b\le\0_m,\quad
    \mu^\star\ge\0_m,\quad
    \mu^\star_i\,(Ay^\star-b)_i=0,\ \ \forall i .
\label{eq:kkt}
\end{equation}
\paragraph{Reference solution and metrics.}
The reference $y^\star$ is computed in two stages. L-BFGS-B first locates the solution and its active set $\mathcal{A}$. The pixels within $\mathcal{A}$ are fixed at their bounds, while
the remaining free pixels are refined by Newton's method.
\begin{table}[H]
\centering
\caption{Selected hyperparameters for image completion using ICNNs.}
\label{tab:hp_icnn}
\begin{tabular}{ll}
\hline
Method & Hyperparameters\\
\hline
COSMIC   & $d t=0.3$, $K=5$, $\beta=0.1$, $\epsilon_\text{sat}=0.1$\\
FL      & $\eta=0.5$, $K=6$\\
Adam-AL & $\alpha_0=10^{-2}$, $\sigma_0=10$, $T=100$, $\gamma_\sigma=1.03$, $\gamma_\alpha=0.97$\\
Adam-PEN & $\alpha_0=10^{-2}$, $\sigma_0=10$, $T=100$, $\gamma_\sigma=1.03$, $\gamma_\alpha=0.97$\\
\hline
\end{tabular}
\end{table}
We present the corresponding convergence of the multiplier to its optimal value in Figure~\ref{fig:icnn_mult}.

We also showcase a sample of the resulting completed faces below.
\begin{figure}[H]\centering\includegraphics[width=0.65\linewidth]{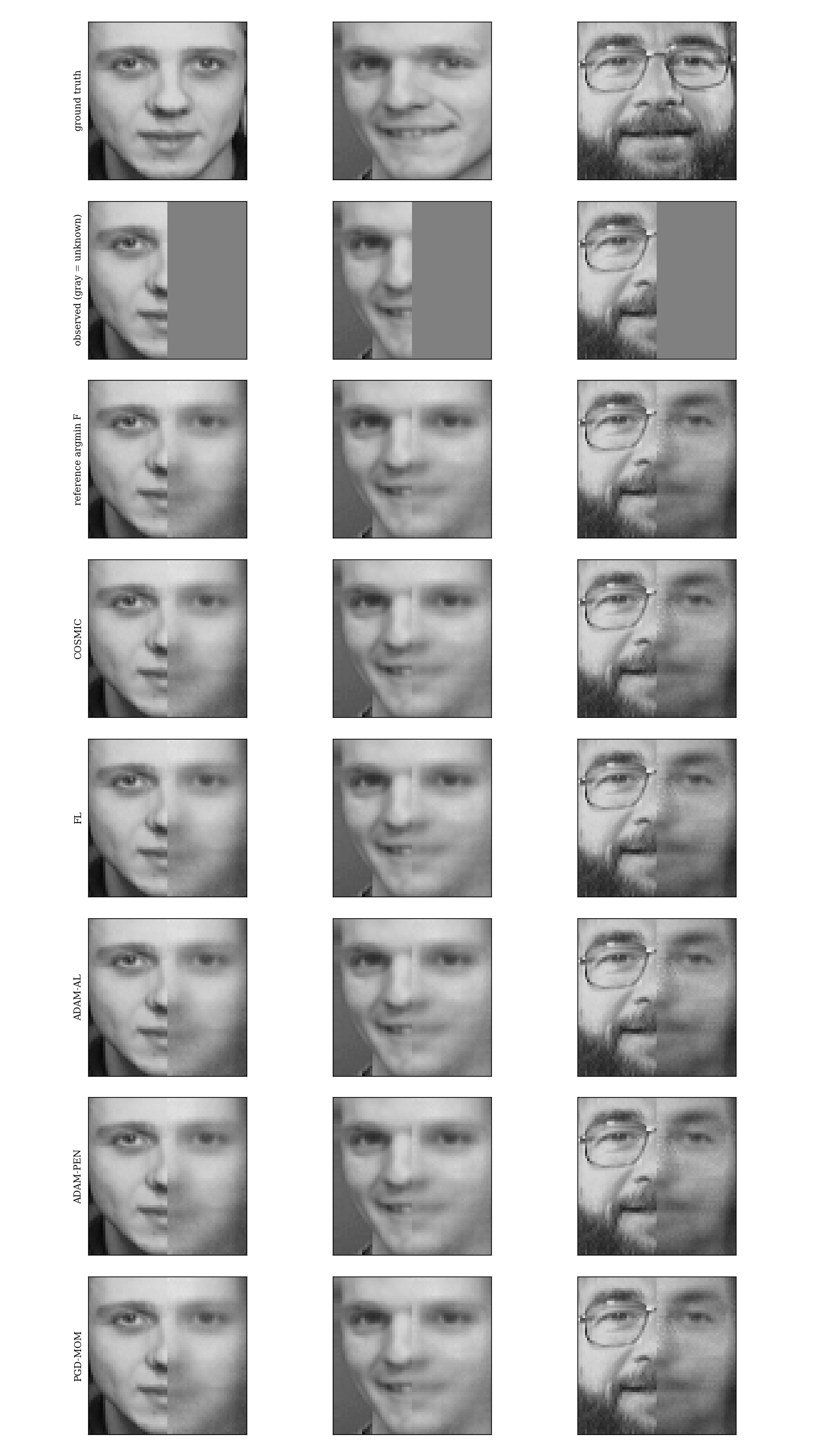}
    \caption{Completed Faces from Olivetti Dataset using \cosmic, FL, Adam-AL, Adam-PEN, and PGD-MOM}
    \label{fig:olivetti_faces}
\end{figure}\label{app:icnn}
\end{document}